\documentclass[12pt,leqno]{amsart}
\usepackage{amssymb,amsmath,amsthm,amscd}
\usepackage{mathrsfs}

\usepackage[all]{xy} 
\CompileMatrices

\numberwithin{equation}{section}

\newcommand{\C}{{\mathbb C}}

\newcommand{\Z}{{\mathbb Z}}
\newcommand{\B}{{\mathscr{B}}}
\newcommand{\I}{{\mathscr{I}}}
\newcommand{\J}{{\mathscr{J}}}
\newcommand{\CM}{{\mathscr{M}}}
\newcommand{\CN}{{\mathscr{N}}}
\newcommand{\CF}{{\mathscr{F}}}
\renewcommand{\L}{{\mathscr{L}}}
\newcommand{\DC}{{\rm D}}
\newcommand{\Dh}[1][X]{{\DC}_\h^\rb(\RD_{#1})}

\newcommand{\RD}{{{\mathscr{D}}}}
\newcommand{\D}{\on{\mathbb{D}}}

\newcommand{\RO}{\mathscr{O}}
\newcommand{\an}{{\mathrm{an}}}
\newcommand{\rb}{{\mathrm{b}}}

\newcommand{\R}{{\rm R}}

\newcommand{\qc}{{\operatorname{qc}}}
\newcommand{\codim}{{\operatorname{codim}}}
\newcommand{\Spec}{{\operatorname{Spec}}}
\newcommand{\Set}{{\mathfrak{S}}}
\newcommand{\CC}{{\mathscr{C}}}
\newcommand{\h}{{\on{hol}}}
\newcommand{\rh}{{\on{rh}}}
\newcommand{\one}{{\bf{1}}}
\newcommand{\seteq}{\operatorname{{:=}}}

\newcommand{\Drh}[1][X]{{\DC}_\rh^{b}(\RD_{#1})}
\newcommand{\qDhge}[1][X]{{}^\Set{\DC}_\h^{\ge0}(\RD_{#1})}
\newcommand{\qDhle}[1][X]{{}^\Set{\DC}_\h^{\le0}(\RD_{#1})}
\newcommand{\qDrhge}[1][X]{{}^\Set{\DC}_\rh^{\ge0}(\RD_{#1})}
\newcommand{\qDrhle}[1][X]{{}^\Set{\DC}_\rh^{\le0}(\RD_{#1})}
\newcommand{\qDge}[1][X]{{}^\Set{\DC}_\qc^{\ge0}(\RD_{#1})}
\newcommand{\qDle}[1][X]{{}^\Set{\DC}_\qc^{\le0}(\RD_{#1})}
\newcommand{\qDh}[2][X]{{}^\Set{\DC}_\hol^{#2}(\RD_{#1})}
\newcommand{\m}{\mathfrak{m}}
\theoremstyle{plain}
\newtheorem{lemma}{Lemma}[section]
\newtheorem{prop}[lemma]{Proposition}
\newtheorem{theorem}[lemma]{Theorem}
\newcommand{\Prop}{\begin{prop}}
\newcommand{\enprop}{\end{prop}}
\newcommand{\Lemma}{\begin{lemma}}
\newcommand{\enlemma}{\end{lemma}}
\newcommand{\Theorem}{\begin{theorem}}
\newcommand{\entheorem}{\end{theorem}}
\newtheorem{corollary}[lemma]{Corollary}
\newcommand{\Cor}{\begin{corollary}}
\newcommand{\encor}{\end{corollary}}
\newtheorem{definition}[lemma]{Definition}

\newcommand{\Def}{\begin{definition}}
\newcommand{\edf}{\end{definition}}

\theoremstyle{definition}
\newtheorem{remark}[lemma]{Remark}
\newtheorem{example}[lemma]{Example}

\newcommand{\Hom}{\operatorname{Hom}}
\newcommand{\End}{\operatorname{End}}

\newcommand{\isoto}[1][]%
{{\mathop{\buildrel{\sim}\over\longrightarrow}\limits_{#1}}}
\renewcommand{\hom}{\operatorname{\it \mathscr{H}\kern-.25em om}}

\newcommand{\Ext}{\operatorname{Ext}}

\newcommand{\Ltens}{\buildrel{\scriptstyle {\mathbb{L}}}\over\otimes}

\newcommand{\tor}{{\on{tor}}}
\newcommand{\M}{{\mathscr M}}
\newcommand{\N}{{\mathscr N}}
\newcommand{\eq}{\begin{eqnarray}}
\newcommand{\eneq}{\end{eqnarray}}

\newcommand{\eqn}{\begin{eqnarray*}}
\newcommand{\eneqn}{\end{eqnarray*}}
\newcommand{\op}{{\operatorname{op}}}
\newenvironment{tenumerate}{
  \begin{enumerate}
  
  }{\end{enumerate}}

\newcommand{\bnum}{\begin{tenumerate}}
\newcommand{\enum}{\end{tenumerate}}
\newcommand{\on}{\operatorname}

\newcommand{\bni}{\begin{tenumerate}}
\newcommand{\eni}{\end{tenumerate}}

\newcommand{\bpf}{\begin{proof}}
\newcommand{\QED}{\end{proof}}
\newcommand{\Proof}{\begin{proof}}
\newcommand{\coh}{{\on{coh}}}

\newcommand{\hol}{{\on{hol}}}

\newcommand{\Supp}{\operatorname{Supp}}
\newcommand{\cl}{\colon}
\newcommand{\To}[1][\phantom{aaaa}]{\xrightarrow{\,#1\,}}

\newcommand{\ba}{\begin{array}}
\newcommand{\ea}{\end{array}}

\renewcommand{\phi}{\varphi}
\newcommand{\bi}{\bni}
\newcommand{\ei}{\eni}

\newcommand{\indlim}[1][]{\varinjlim\limits_{#1}}
\newcommand{\set}[2]{\left\{#1\,;\,#2\right\}}
\newcommand{\supp}{\operatorname{supp}}
\newcommand{\Mod}{\operatorname{Mod}}
\newcommand{\hs}{\hspace}

\newcommand{\eqsub}{\begin{subequations}\begin{eqnarray}}
\newcommand{\eneqsub}{\end{eqnarray}\end{subequations}}

\newcommand{\qDA}[2][]{{}^{\Phi{#1}}\DC^{#2}_\qc(\A{#1})}
\newcommand{\sDA}[2][]{{}^\Set\DC^{#2}_\qc(\A{#1})}
\newcommand{\qA}[1]{{}^\Phi\DC^{#1}(\A)}

\newcommand{\ol}{\overline}
\newcommand{\Tri}[1][]{{{\mathbf{D}}^{#1}}}
\newcommand{\A}{\mathscr{A}}
\newcommand{\Sol}{{\on{Sol}}}
\newcommand{\sDO}[2][X]{{}^\Set\DC^{#2}_\qc(\RO_{#1})}
\newcommand{\qDO}[2][X]{{}^\Phi\DC^{#2}_\qc(\RO_{#1})}
\newcommand{\qDOX}[2][\Phi]{{}^{#1}\DC^{#2}_\qc(\RO_X)}
\newcommand{\sDOc}[2][X]{{}^\Set\DC^{#2}_\coh(\RO_{#1})}
\newcommand{\qDOc}[2][\Phi]{{}^{#1}\DC^{#2}_\coh(\RO_X)}
\newcommand{\DOc}[1][\rb]{\DC^{#1}_\coh(\RO_{X})}
\newcommand{\DO}[1][\rb]{\DC^{#1}_\qc(\RO_{X})}
\newcommand{\sDD}[1][\rb]{{}^\Set\DC^{#1}_\qc(\RD_{X})}
\renewcommand{\le}{\leqslant}
\renewcommand{\ge}{\geqslant}
\newcommand{\Ind}{{{\on{Ind}}}}
\newcommand{\Pro}{{{\on{Pro}}}}
\newcommand{\ko}{{{\mathbf{k}}}}
\newcommand{\ext}{{\mathscr{E}xt}}

\newcommand{\sd}{support datum}
\newcommand{\sds}{support data}
\newcommand{\tv}{\mathbb{T}}
\newcommand{\SH}{\mathscr{H}}

\begin{document}

\title[t-structures on the derived categories of D-modules and O-modules]%
{t-structures on the derived categories of holonomic D-modules
and coherent O-modules}
\author{Masaki KASHIWARA}
\address{Research Institute for Mathematical Sciences,
Kyoto University, Kyoto 606, Japan
}
\thanks{This research is partially supported by 
Grant-in-Aid for Scientific Research (B1)13440006,
Japan Society for the Promotion of Science.}
\dedicatory{Dedicated to Boris Feigin on his fiftieth birthday}

\begin{abstract}
We give the description of the $t$-structure 
on the derived category of regular holonomic $\RD$-modules
corresponding to the trivial $t$-structure on
the derived category of constructible sheaves
via Riemann-Hilbert correspondence.
We give also the condition for a decreasing sequence
of families of supports to give a $t$-structure on the derived category of 
coherent $\RO$-modules.
\end{abstract}

\maketitle

\section{introduction}
It was one of the motivations of the introduction
of the notions of $t$-structures and perverse sheaves
by A. Beilinson, J. Bernstein, P. Deligne and O. Gabber
(\cite{FP})
to ask what are the objects corresponding 
to regular holonomic $\RD_X$-modules
by the Riemann-Hilbert correspondence (\cite{K})
$\smash{\R\hom_{\RD_X}(-,\RO_X)\cl\DC^\rb_\rh(\RD_X)\isoto
\DC^\rb_c(\C_X)^\op}$.
Here $X$ is a complex manifold, $\DC^\rb_c(\C_X)^\op$
is the opposite category of the derived category $\DC^\rb_c(\C_X)$
of bounded complexes of
sheaves on $X$ with constructible cohomologies,
and $\DC^\rb_\rh(\RD_X)$ is the derived category of 
bounded complexes of $\RD_X$-modules with regular holonomic cohomologies.
With the notion of $t$-structures,
the answer is: the $t$-structure of the middle perversity
on $\DC^\rb_c(\C_X)$ corresponds to 
the trivial $t$-structure on $\DC^\rb_\rh(\RD_X)$.

The purpose of this paper is to answer the converse question:
what is the $t$-structure on
$\DC^\rb_\rh(\RD_X)$ corresponding to the trivial $t$-structure on
$\DC^\rb_c(\C_X)$.
In fact, this $t$-structure can be extended to a $t$-structure
on the derived category $\DC^\rb_\qc(\RD_X)$
of bounded complexes of $\RD_X$-modules with quasi-coherent cohomologies.
In this paper, 
we treat an algebraic case, i.e.\ when $X$ is a smooth algebraic variety.
In this case, the corresponding
$t$-structure $(\sDD[\le0],\sDD[\ge0])$ is given as follows.
\eqn
&&\sDD[\le0]=\set{\CM\in\DC^\rb_\qc(\RD_X)}%
{\mbox{$\codim \Supp(H^n(\CM))\ge n$ for every $n\ge0$}}\,,\\
&&\sDD[\ge0]=\set{\CM\in\DC^\rb_\qc(\RD_X)}%
{\parbox{210pt}{$\SH^k_Z(\CM)=0$ for any closed subset $Z$ and $k<\codim Z$}}.
\eneqn

Similar results hold for a complex analytic manifold $X$.

More general results are given in this paper.
For a left Noetherian $\RO_X$-ring $\A$ quasi-coherent over $\RO_X$
and a decreasing sequence of families of supports $\Phi=\{\Phi^n\}_{n\in\Z}$,
the pair
{\allowdisplaybreaks
\eqn
&&{}^\Phi\DC^{\le0}_\qc(\A)\seteq
\set{M\in\DC^\rb_\qc(\A)}{\mbox{$\Supp(H^k(M))\subset\Phi^{k}$ for any $k$}},
\\[5pt]
&&{}^\Phi\DC^{\ge0}_\qc(\A)\seteq
\set{M\in\DC^\rb_\qc(\A)}{\mbox{$\SH_{\Phi^n}^k(M)=0$ if $k<n$}}
\eneqn
}%
gives a $t$-structure on $\DC^\rb_\qc(\A)$ (Theorem \ref{th:t-str}).
This construction is similar to the one of the perverse
sheaves.
The $t$-structure $(\sDD[\le0],\sDD[\ge0])$ corresponds to the case where
$\A=\RD_X$ and $\Phi=\Set$,
where $\Set^k$ is the family of supports consisting of 
closed subsets with codimension $\ge k$.

However, this $t$-structure does not
always induce a $t$-structure on $\DC^\rb_\coh(\RO_X)$, the derived category
of complexes of $\RO_X$-modules with coherent cohomologies.
We give the
necessary and sufficient condition for $\Phi$
to give a $t$-structure on $\DC^\rb_\coh(\RO_X)$ (Theorem \ref{th:tOc}).
\footnote%
{\normalsize After writing this paper, the author was informed the existence of the paper
{\em ``Perverse coherent sheaves {\rm(}after Deligne{\rm)''}}, math.AG/0005152 by Roman Bezrukavnikov,
in which it is proved in a more general setting that the condition
for $\Phi$
to give a $t$-structure on $\DC^\rb_\coh(\RO_X)$ is sufficient.}
This condition resembles the one for perversity.

\smallskip
The paper is organized as follows.
In \S \ref{sect: tstr},
we briefly review the notion of $t$-structures.

In \S \ref{sect:noet}, we introduce 
the $t$-structure $({}^\Phi\DC^{\le0}_\qc(\A),{}^\Phi\DC^{\ge0}_\qc(\A))$
on the triangulated category $\DC^\rb_\qc(\A)$,
and studies
its properties.

In \S \ref{sect:dual}, we study the
$t$-structure on the derived category $\DC^\rb_\coh(\RO_X)$
which corresponds to the standard $t$-structure on $\DC^\rb_\coh(\RO_X)$
by the duality functor $\R\hom_{\RO_X}(-,\RO_X)$.
We also give the relation between this $t$-structure and flatness
(Proposition \ref{prop:flat}).

In \S \ref{sect:tsto}, we give the
condition for a decreasing sequence of families of supports
to give a $t$-structure on $\DC^\rb_\coh(\RO_X)$ (Theorem \ref{th:tOc}).
This is a generalization of \cite[Exercise X.2]{KS}.

In \S \ref{sec:D}, we study the $t$-structure on
$\DC^\rb_\qc(\RD_X)$.

In the last section, we give a proof of the stability of injectivity 
under filtrant inductive limits (Lemma \ref{lem:indinj}).

\medskip
While writing this paper, the author received the
preprint \cite{Y} by A. Yekutieli and James J. Zhang, 
where they proposed similar $t$-structures.

\section{t-structure}\label{sect: tstr}
In this section, we recall the notion of $t$-structures.
For details, see \cite{FP, KS}.
Let $\Tri$ be a triangulated category.
Let $\Tri[\le0]$ and $\Tri[\ge0]$ be full subcategories of $\Tri$.
We set $\Tri[\le n]\seteq\Tri[\le0][-n]$
and $\Tri[\ge n]\seteq\Tri[\ge0][-n]$.
We also use the notations:
$\Tri[<n]\seteq\Tri[\le n-1]$ and 
$\Tri[>n]\seteq\Tri[\ge n+1]$.
The pair $(\Tri[\le0],\Tri[\ge0])$ is called
a {\em $t$-structure} on $\Tri$ if it satisfies the following conditions:
\begin{subequations}
\eq
&&\mbox{$\Tri[<0]\subset\Tri[\le0]$ and $\Tri[>0]\subset\Tri[\ge0]$,}
\label{eq:1a}
\\[5pt]
&&\mbox{$\Hom_{\Tri}(X,Y)=0$ for $X\in\Tri[<0]$ and $Y\in\Tri[\ge0]$,}
\label{eq:1b}\\[5pt]
&&
\ba{l}
\parbox{370pt}{For any $X\in\Tri$, there exists a distinguished triangle
$X'\to X\to X''\To[+1]$
with $X'\in\Tri[<0]$ and $X''\in\Tri[\ge0]$.}
\ea\label{eq:1c}
\eneq
Then one has
\eq\label{eq:1d}
&&\parbox{350pt}{For a distinguished triangle
$X'\to X\to X''\To[+1]$, if $X'$ and $X''$ belong to
$\Tri[\le0]$ (resp.\ $\Tri[\ge0]$), then so does $X$.}
\eneq
\end{subequations}

Note that the $t$-structure is a self-dual notion: 
if $(\Tri[\le0],\Tri[\ge0])$
is a $t$-structure on a triangulated category $\Tri$,
then $((\Tri[\ge0])^\op,(\Tri[\le0])^\op)$
is a $t$-structure on the opposite triangulated category $\Tri^\op$.

In the paper, we use the following results.

\Lemma\label{lem:t-str}\label{lem:t-strd}
Let $\Tri$ be a triangulated category,
and assume that $(\Tri[\le0],\Tri[\ge0])$
satisfies \eqref{eq:1d}.
Let $M_1\to M\to M_2\To[+1]$ be a distinguished triangle.
If one of the following conditions {\rm (i)} and {\rm (ii)} is satisfied,
then there exists a distinguished triangle
$M'\to M\to M''\To[+1]$
with $M'\in \Tri[<0]$ and $M''\in\Tri[\ge0]$.
\bi
\item
$M_1\in \Tri[<0]$ and there exists a distinguished triangle
$M'_2\to M_2\to M''_2\To[+1]$ with $M'_2\in \Tri[<0]$
and $M''_2\in \Tri[\ge0]$.
\item
$M_2\in \Tri[\ge0]$ and there exists a distinguished triangle
$M'_1\to M_1\to M''_1\To[+1]$ with $M'_1\in \Tri[<0]$
and $M''_1\in \Tri[\ge0]$.
\ei
\enlemma
\proof
Assume the condition (i).
By the octahedral axiom,
$$
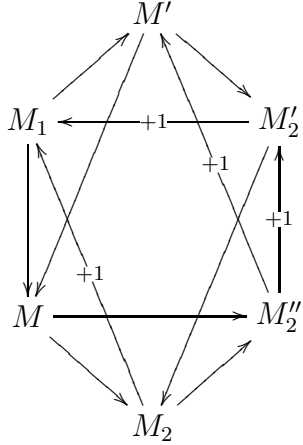
\xymatrix{
&M'\ar[dr]\ar[dddl]\\
M_1\ar[ur]\ar[dd]&&M_2'\ar[ll]|{+1}\ar[dddl]\\\\
M\ar[dr]\ar[rr]&&M''_2\ar[uu]|{+1}\ar[uuul]|{+1}\\
&M_2\ar[ur]\ar[uuul]|{+1}
}$$
one has distinguished triangles
$M_1\to M'\to M'_2\To[+1]$ and $M'\to M\to M''_2\To[+1]$.
Hence $M'\in \Tri[<0]$ by \eqref{eq:1d}.

The condition (ii) is the dual form of (i).
\qed

\medskip
The following lemma is almost obvious.
\Lemma\label{lem:equi}
Let $F\cl \Tri\to\Tri'$ be an equivalence of triangulated categories.
Let $(\Tri[\le0],\Tri[\ge0])$ be a $t$-structure on $\Tri$,
and let $(\Tri'{}^{\le0},\Tri'{}^{\ge0})$ be a pair of full subcategories of
$\Tri'$
satisfying \eqref{eq:1a} and \eqref{eq:1b}.
If $F$ sends $\Tri[\le0]$ to $\Tri'{}^{\le0}$
and $\Tri[\ge0]$ to $\Tri'{}^{\ge0}$,
then $(\Tri'{}^{\le0},\Tri'{}^{\ge0})$ is a $t$-structure on $\Tri'$.
\enlemma

Let $\Tri'$ be another triangulated category,
and let $(\Tri[\le0],\Tri[\ge0])$
and $(\Tri'{}^{\le0},\Tri'{}^{\ge0})$ be a $t$-structure on $\Tri$ and
$\Tri'$, respectively.
A functor $F\cl \Tri\to \Tri'$ of
triangulated categories
is called {\em left exact} if $F$ sends
$\Tri[\ge0]$ to $\Tri'{}^{\ge0}$.
It is called {\em right exact} if 
it sends $\Tri[\le0]$ to $\Tri'{}^{\le0}$.
It is called {\em exact} if it is left exact and right exact.

\section{t-structures on a Noetherian scheme}\label{sect:noet}
Let $X$ be a finite-dimensional Noetherian separated scheme.
Let $\A$ be an $\RO_X$-ring, i.e.\ 
a (not necessarily commutative) ring
on $X$ with a ring morphism $\RO_X\to\A$.
We assume that $\A$ is quasi-coherent as a left $\RO_X$-module,
and that  $\A$ is left Noetherian (e.g.\ see 
\cite[Definition A.7]{K1}). Under the first assumption, the second
is equivalent to saying that $\A(U)$ is a left
Noetherian ring for any affine open subset $U$.

Let $\Mod(\A)$ be the category of (left) $\A$-modules, and
$\Mod_\qc(\A)$ (resp.\ $\Mod_\coh(\A)$) 
its full subcategory of quasi-coherent (resp.\ coherent) $\A$-modules.
Note that an $\A$-module is quasi-coherent over $\A$ 
if and only if it is quasi-coherent over $\RO_X$.

One has the following lemma whose proof is given in 
the last section \ref{sec:proof}.

\Lemma\label{lem:indinj}
Any filtrant inductive limit of injective objects of
$\Mod(\A)$ is injective.
\enlemma

Let $\DC(\A)$ be the derived category of
$\Mod(\A)$, and
let $\DC^\rb(\A)$ be the full subcategory of
$\DC(\A)$ consisting of objects with bounded cohomologies.
Let $\DC^\rb_\qc(\A)$ (resp.\ $\DC^\rb_\coh(\A)$) be the full subcategory of
$\DC^\rb(\A)$ consisting of objects
with quasi-coherent (resp.\ coherent) cohomologies.
We define
\eqn
&&\DC^{\le n}(\A)\seteq\set{M\in\DC^\rb(\A)}{\mbox{$H^k(M)=0$ for $k>n$}},\\
&&\DC^{\ge n}(\A)\seteq\set{M\in\DC^\rb(\A)}{\mbox{$H^k(M)=0$ for $k<n$}}
\eneqn
and similarly $\DC^{\le n}_\qc(\A)$, $\DC^{\le n}_\coh(\A)$, etc.
We call the $t$-structure $(\DC^{\le0}(\A),\DC^{\ge0}(\A))$
the {\em standard $t$-structure}.
Let us denote by
$\tau^{\le n}\cl\DC^\rb(\A)\to\DC^{\le n}(\A)$ and
$\tau^{\ge n}\cl\DC^\rb(\A)\to\DC^{\ge n}(\A)$
the truncation functors with respect to the  standard $t$-structure.

In this section, we shall give $t$-structures on
$\DC^\rb(\A)$ and $\DC^\rb_\qc(\A)$.

\medskip
In this paper, a {\em family of supports} 
means a set $\Phi$
of closed subsets of $X$ satisfying the following conditions:
\bi
\item If $Z\in \Phi$ and $Z'$ is a closed subset of $Z$, 
then $Z'\in\Phi$.\label{it}
\item For $Z$, $Z'\in \Phi$, the union
$Z\cup Z'$ belongs to $\Phi$.
\item $\Phi$ contains the empty set.
\ei

Then the set of families of supports has the structure of an ordered set
by inclusion.
The join of families of supports $\Phi_1$ and $\Phi_2$
is as follows:
\eqn
\Phi_1\cup\Phi_2&=&\set{Z}{\parbox{230pt}{$Z$ is a closed subset
such that $Z\subset Z_1\cup Z_2$ 
for some $Z_1\in\Phi_1$ and some $Z_2\in\Phi_2$}}\\
&=&\set{Z}{\parbox{230pt}{$Z$ is a closed subset
such that $Z=Z_1\cup Z_2$ 
for some $Z_1\in\Phi_1$ and some $Z_2\in\Phi_2$}}.
\eneqn
Note that an irreducible closed subset $Z$ is
a member of $\Phi_1\cup\Phi_2$ if and only if
one has either
$Z\in \Phi_1$ or $Z\in \Phi_2$.
Their meet is given by
\eqn
&&\Phi_1\cap\Phi_2=\set{Z}{\mbox{$Z\in\Phi_1$ and $Z\in\Phi_2$}}.
\eneqn
We sometimes regard a closed subset $S$ of $X$ as the family of supports
consisting of closed subsets of $S$.
Hence for a closed subset $S$ and a family of supports
$\Phi$, one has
$\Phi\cap S=\set{Z\in\Phi}{Z\subset S}$, and
$\Phi\cup S$ is the the family of supports
consisting of closed subsets $Z$ such that
$\ol{Z\setminus S}\in\Phi$.


\medskip
For a sheaf $F$, one sets
$\Gamma_\Phi(F)=\indlim[Z\in\Phi]\Gamma_Z(F)$.
Then one has
$$\Gamma(U;\Gamma_\Phi(F))=\set{s\in F(U)}{\overline{\supp(s)}\in\Phi}
\mbox{ for any open subset $U$,}$$
because $X$ is Noetherian.
Then $\Gamma_\Phi$ is a left exact functor from
$\Mod(\A)$ to itself, and also it sends $\Mod_\qc(\A)$ to itself.
It commutes with filtrant inductive limits.
It is well-known that
$\Gamma_Z$ sends injective objects of $\Mod(\A)$ to injective objects.
Hence, by Lemma \ref{lem:indinj},
the functor
$\Gamma_\Phi$ also sends injective objects of $\Mod(\A)$ to injective objects.
Let $\R\Gamma_\Phi\cl \DC^\rb(\A)\to \DC^\rb(\A)$
be the right derived functor of $\Gamma_\Phi$.
Similarly, we can define $\R\Gamma_\Phi\cl\DC^\rb(\Z_X)\to\DC^\rb(\Z_X)$,
where $\DC^\rb(\Z_X)$ is the derived category of bounded 
complexes of $\Z_X$-modules. The functors $\R\Gamma_\Phi$ commutes 
with the forgetful functor $\DC^\rb(\A)\to \DC^\rb(\Z_X)$.

For $M\in\DC^\rb(\A)$, let $\Supp(M)$ denote the family of supports 
consisting of closed subsets $Z$ such that
$Z\subset\bigcup_{i=1}^n\overline{\supp(s_i)}$ 
for open subsets $U_i$, integers $n_i$ and
$s_i\in H^{n_i}(M)(U_i)$.
Hence, $\Supp(M)\subset \Phi$ is equivalent to saying that
$\R\Gamma_\Phi(M)\isoto M$.

One sets
$\SH^n_\Phi(M)\seteq H^n(\R\Gamma_\Phi(M))$.
One has $\SH^n_\Phi(M)\simeq\smash{\indlim[Z\in\Phi]}\SH^n_Z(M)$
for any $M\in \DC^\rb(\A)$ and any integer $n$.

For any $M\in \DC^\rb_\qc(\A)$, one has an isomorphism (\cite{G})
\eq
&&\SH^n_\Phi(M)\simeq\indlim[I]
\ext^n_{\RO_X}(\RO_X/I,M)\eneq
as $\RO_X$-modules.
Here the inductive limit is taken over
coherent ideals $I$ of $\RO_X$
such that $\Supp(\RO_X/I)\in\Phi$.
Since an $\A$-module is quasi-coherent if and only if 
it is quasi-coherent over $\RO_X$,
the functor $\R\Gamma_\Phi$  induces a functor
$$\R\Gamma_\Phi\cl \DC^\rb_\qc(\A)\to \DC^\rb_\qc(\A).$$

\begin{remark}\label{rem:flabby}
Although one neither gives a proof nor uses it in this paper,
one has
$\SH^n_\Phi(M)\simeq\indlim[I]
\ext^n_{\A}(\A/\A I,M)$.
Here the inductive limit is taken over
$I$ as above.
Hence an injective object of 
$\Mod_\qc(\A)$ is a flabby sheaf (see Remark \ref{rem:injinj}).
\end{remark}

\medskip

The following lemma is obvious.
\Lemma\label{lem:aux}
\bi
\item For families of supports $\Phi$ and $\Phi'$, one has
$\R\Gamma_\Phi\circ\R\Gamma_{\Phi'}=\R\Gamma_{\Phi\cap\Phi'}$.
\item
For $M$, $K\in\DC^\rb(\A)$ with $\Supp(M)\subset\Phi$, 
one has
$$\Hom_{\DC^\rb(\A)}(M,\R\Gamma_\Phi K)\isoto
\Hom_{\DC^\rb(\A)}(M,K)$$
and $\ext^n_\A(M,\R\Gamma_\Phi K)\isoto\ext^n_\A(M,K)$ for every $n$.
\item
If $M\in\DC^{\ge n}(\A)$,
then
$\R\Gamma_\Phi(M)\in\DC^{\ge n}(\A)$
and $H^n(\R\Gamma_\Phi(M))=\Gamma_\Phi(H^n(M))$.
\item \label{en:3}
For a locally closed subset $Z$ of $X$ and a family $\Phi$
of supports, one has
$$\R\Gamma_Z\circ\R\Gamma_\Phi\simeq \R\Gamma_\Phi\circ\R\Gamma_Z.$$
\item
For an open embedding $j\cl U\hookrightarrow X$, one has
$\R j_*\circ\R\Gamma_{\Phi|_U}\simeq \R\Gamma_\Phi\circ\R j_*$.
Here $\Phi|_U$ is the family of supports on $U$
given by $$\qquad\Phi|_U\seteq
\set{Z}{\mbox{$Z$ is a closed subset of $U$ such that
$\ol{Z}\in\Phi$}}=\set{Z\cap U}{Z\in\Phi}.$$
\ei
\enlemma

\Lemma
Let $M\in \DC^\rb(\A)$ and $n$ an integer, and
let $\Phi$ and $\Psi$ be families of supports.
\bi
\item
If $\Psi\subset\Phi$ and
$\R\Gamma_\Phi(M)\in \DC^{\ge n}(\A)$,
then $\R\Gamma_\Psi(M)\in \DC^{\ge n}(\A)$.
\item
$\R\Gamma_\Phi(M)\in \DC^{\ge n}(\A)$
if and only if
$\R\Gamma_Z(M)\in \DC^{\ge n}(\A)$ for every $Z\in\Phi$.
\item
If $\R\Gamma_\Phi(M)$, $\R\Gamma_\Psi(M)\in \DC^{\ge n}(\A)$,
then
$\R\Gamma_{\Phi\cup\Psi}(M)\in \DC^{\ge n}(\A)$.
\item
If $\SH^k_{{\ol{\{{x}\}}}}(M)_x=0$ for any $k<n$ and
any $x$ such that $\ol{\{x\}}\in\Phi$,
then
$\R\Gamma_\Phi(M)\in \DC^{\ge n}(\A)$.
\ei
\enlemma

\proof\quad
(i) follows immediately from
$\R\Gamma_\Psi(M)\simeq\R\Gamma_\Psi\R\Gamma_\Phi(M)$, and
%
(ii) follows from (i) and
$\SH^k_\Phi(M)\simeq\indlim[Z\in\Phi]\SH^k_Z(M)$.

Let us show (iii).
By (i) and (ii), it is enough to show that
$\R\Gamma_{Z_1\cup Z_2}(M)\in \DC^{\ge n}(\A)$
for closed subsets
$Z_1$ and $Z_2$ such that
$\R\Gamma_{Z_1}(M)$, $\R\Gamma_{Z_2}(M)\in \DC^{\ge n}(\A)$.
Since one has
$$\tau^{<n}\R\Gamma_{Z_1\cup Z_2}(M)|_{X\setminus Z_1}
\simeq\tau^{<n}\R\Gamma_{Z_2}(M)|_{X\setminus Z_1}=0,$$
the support of $\tau^{<n}\R\Gamma_{Z_1\cup Z_2}(M))$ is contained in $Z_1$.
In the distinguished triangle
$$\R\Gamma_{Z_1}\bigl(\tau^{<n}\R\Gamma_{Z_1\cup Z_2}(M)\bigr)\to
\R\Gamma_{Z_1}
\R\Gamma_{Z_1\cup Z_2}(M)\to\R\Gamma_{Z_1}
\bigl(\tau^{\ge n}\R\Gamma_{Z_1\cup Z_2}(M)\bigr),$$
the second term is isomorphic to
$\R\Gamma_{Z_1}(M)$ which belongs to $\DC^{\ge n}(\A)$,
and the last term also belongs to
$\DC^{\ge n}(\A)$.
Hence the first term, isomorphic to $\tau^{<n}\R\Gamma_{Z_1\cup Z_2}(M)$,
belongs to $\DC^{\ge n}(\A)$, and therefore it vanishes.

\smallskip\noindent
(iv)\quad By the induction on $n$, one may assume that
$\R\Gamma_\Phi(M)\in \DC^{\ge n-1}(\A)$.
For an open set $U$ and $s\in\Gamma(U;\SH^{n-1}_\Phi(M))$,
set $S\seteq\supp(s)$.
Let us show that $S=\emptyset$.
Otherwise, let us take the generic point $x$ 
of an irreducible component of $S$.
Then one has $\ol{\{x\}}\in\Phi$.
Since $\Gamma_S\bigl(\SH^{n-1}_\Phi(M)\bigr)=\SH^{n-1}_S(M)$
and $\SH^{n-1}_S(M)_x=\SH^{n-1}_{\ol{\{x\}}}(M)_x$,
the germ of $s$ at $x$ must vanish, which is a contradiction.
\qed

\medskip
A {\em \sd}\ 
is a decreasing sequence $\Phi\seteq\{\Phi^n\}_{n\in \Z}$ of 
families of supports satisfying the following conditions:
\begin{subequations}
\eq
&&\mbox{for $n\ll0$, $\Phi^n$ is the set of all closed subsets of $X$,}\\
&&\mbox{for $n\gg0$, $\Phi^n$ is $\{\emptyset\}$.}
\eneq\end{subequations}
For a \sd\ $\Phi$ and an integer $n$, let $\sigma^{\le n}\Phi$ denote
the \sd\  given by
\eqn&&
(\sigma^{\le n}\Phi)^k=
\left\{
\ba{ll}
\Phi^k&\mbox{for $k\le n$},\\
\{\emptyset\}&\mbox{for $k>n$}.
\ea\right.
\eneqn
Let $\tv$ denote the \sd\ given by
$$\tv^n=
\left\{
\ba{ll}
\left\{\mbox{all closed subsets of $X$}\right\}&\mbox{for $n\le0$,}\\
\{\emptyset\}&\mbox{for $n>0$.}
\ea
\right.$$

For a \sd\ $\Phi=\{\Phi^n\}_{n\in \Z}$, we set
\eq
&&\begin{array}{l}
{}^\Phi\DC^{\le n}(\A)\seteq
\set{M\in\DC^\rb(\A)}{\mbox{$\Supp(H^k(M))\subset\Phi^{k-n}$ for any $k$}},
\\[5pt]
{}^\Phi\DC^{\ge n}(\A)\seteq
\set{M\in\DC^\rb_\qc(\A)}{\mbox{$\R\Gamma_{\Phi^k}(M)\in
\DC^{\ge k+n}(\A)$ for any $k$}},
\ea
\eneq
and $\qDA{\le n}\seteq{}^\Phi\DC^{\le n}(\A)\cap\DC^\rb_\qc(\A)$,
$\qDA{\ge n}\seteq{}^\Phi\DC^{\ge n}(\A)\cap\DC^\rb_\qc(\A)$.

\Theorem\label{th:t-str}
${}^\Phi\DC^{\rb}(\A)\seteq
({}^\Phi\DC^{\le 0}(\A),{}^\Phi\DC^{\ge0}(\A))$
and
$\qDA{\rb}\seteq
(\qDA{\le0},\qDA{\ge0})$
are a $t$-structure on $\DC^\rb(\A)$ and $\DC^\rb_\qc(\A)$, respectively.
\entheorem

Note that ${}^\tv\DC^\rb(\A)$ is the standard $t$-structure.

\medskip
Since the proof for ${}^\Phi\DC^{\rb}_\qc(\A)$ is similar, we only
give a proof for
$\qA{\rb}$.
We divide the proof of the theorem into several steps.
It is evident that $\qA{\rb}$
satisfies the conditions \eqref{eq:1a} and \eqref{eq:1d}.
Let us show that it satisfies \eqref{eq:1b}.
\Lemma
For $M\in\qA{<0}$ and
$K\in\qA{\ge0}$, one has
\eqn\label{eq:fn}
\Hom_{\DC^\rb(\A)}(M,K)=0\mbox{ and $\ext^n_\A(M,K)=0$ for $n\le0$.}
\eneqn
\enlemma
\proof
We shall show
\begin{equation}
\mbox{$\Hom_{\DC^\rb(\A)}(M,K)=0\,$
for $M\in \qA{<0}\cap \DC^{\le n}(\A)$
and $K\in \qA{\ge0}$}\tag*{(\ref{eq:fn})${}_n$}
\end{equation}
by the induction on $n$.
Since $\qA{\ge0}\subset\DC^{>n}(\A) $
for $n\ll 0$, (\ref{eq:fn})${}_n$ is satisfied for $n\ll0$.
Assuming (\ref{eq:fn})${}_{n-1}$, let us show
(\ref{eq:fn})${}_n$.
One has
$\Supp(H^n(M)[-n])\subset \Phi^{n+1}$, 
$H^n(M)[-n]\in \DC^{\le n}(\A)$
and $\R\Gamma_{\Phi^{n+1}}K\in \DC^{>n}(\A)$.
Hence one has
$$\Hom_{\DC^\rb(\A)}(H^n(M)[-n],K)
\simeq\Hom_{\DC^\rb(\A)}(H^n(M)[-n],\R\Gamma_{\Phi^{n+1}}(K))=0.$$
The distinguished triangle
$\tau^{<n}M\to M\to H^n(M)[-n]\To[+1]$
induces an exact sequence
$$\Hom_{\DC^\rb(\A)}(H^n(M)[-n],K)\to
\Hom_{\DC^\rb(\A)}(M,K)\to\Hom_{\DC^\rb(\A)}(\tau^{<n}M,K).$$
Since $\Hom_{\DC^\rb(\A)}(\tau^{<n}M,K)=0$ by (\ref{eq:fn})${}_{n-1}$,
we obtain $\Hom_{\DC^\rb(\A)}(M,K)=0$.

The proof of the second statement is similar.
\qed

\medskip
Hence, it remains to prove
the following statement
for any $M\in \DC^\rb(\A)$:
\eq\label{eq:dist}
&&\parbox{400pt}{There exists a distinguished triangle
$M'\to M\to M''\To[+1]$
with $M'\in\qA{<0}$
and $M''\in\qA{\ge0}$.}
\eneq

\proof[Proof of \eqref{eq:dist}]
For $n\in \Z$, let us consider the following statement:
\begin{equation}
\parbox{400pt}{
The property \eqref{eq:dist} holds
if $M\in \DC^\rb(\A)$ satisfies
$\R\Gamma_{\Phi^i}(M)\in \DC^{\ge i}(\A)$ for any $i\le n$,
i.e.\ if $M\in{}^{\sigma^{\le n}\Phi}\DC^{\ge0}(\A)$.}
\tag*{(\ref{eq:dist})${}_n$}
\end{equation}
Since every $M\in \DC^\rb(\A)$
satisfies $\R\Gamma_{\Phi^i}(M)\in \DC^{\ge i}(\A)$
for $i\ll0$, it is enough to show (\ref{eq:dist})${}_n$
for every $n$.
We shall show (\ref{eq:dist})${}_n$ by the descending induction on $n$.
Note that \eqref{eq:dist}${}_n$
holds for sufficiently large $n$ such that $\Phi^{n+1}=\{\emptyset\}$,
because such an $M$ satisfies $M\in\qA{\ge0}$.

We shall prove \eqref{eq:dist}${}_{n}$ 
by assuming \eqref{eq:dist}${}_{n+1}$.
Suppose that $M\in \DC^\rb(\A)$ satisfies
$\R\Gamma_{\Phi^i}(M)\in \DC^{\ge i}(\A)$ for any $i\le n$.
Let us consider a distinguished triangle
$$\tau^{\le n}\R\Gamma_{\Phi^{n+1}}(M)\to M\to M''\To[+1].$$
Since $\R\Gamma_{\Phi^{n+1}}(M)=\R\Gamma_{\Phi^{n+1}}\R\Gamma_{\Phi^{n}}(M)$,
one has
$\R\Gamma_{\Phi^{n+1}}(M)\in \DC^{\ge n}(\A)$
and
\eq\label{eq:k=n}
&&\SH^n_{\Phi^{n+1}}(M)\simeq
\Gamma_{\Phi^{n+1}}(\SH^n_{\Phi^{n}}(M)).
\eneq
Moreover one has
\eqn
\Supp(H^i(\tau^{\le n}\R\Gamma_{\Phi^{n+1}}(M)))
&=&\left\{
\ba{ll}
\emptyset&\mbox{for $i\not=n$,}\\
\subset\Phi^{n+1}&\mbox{for $i=n$,}
\ea
\right.
\eneqn
One concludes therefore that
$\tau^{\le n}\R\Gamma_{\Phi^{n+1}}(M)\in\qA{<0}$.
If one shows
\eq\label{eq:M''}
\mbox{$\R\Gamma_{\Phi^i}(M'')\in \DC^{\ge i}(\A)$ for any $i\le n+1$,}
\eneq
then we can apply \eqref{eq:dist}${}_{n+1}$ to $M''$ and 
\eqref{eq:dist} is satisfied for $M''$.
Hence Lemma \ref{lem:t-str} implies that
\eqref{eq:dist} is satisfied for $M$.

It remains to prove \eqref{eq:M''}, i.e.\ $H^k(\R\Gamma_{\Phi^i}(M''))=0$
for $k<i\le n+1$.
For $i\le n+1$ one has
$\R\Gamma_{\Phi^i}(\tau^{\le n}\R\Gamma_{\Phi^{n+1}}(M))
=\tau^{\le n}\R\Gamma_{\Phi^{n+1}}(M)$.
Hence we have a long exact sequence
\eqn
&&H^k(\tau^{\le n}\R\Gamma_{\Phi^{n+1}}(M))
\To[\xi] H^k(\R\Gamma_{\Phi^i}(M))\to H^k(\R\Gamma_{\Phi^i}(M''))\to\\
&&\hs{100pt}
\to H^{k+1}(\tau^{\le n}\R\Gamma_{\Phi^{n+1}}(M))
\To[\eta] H^{k+1}(\R\Gamma_{\Phi^i}(M))\,.
\eneqn
Hence it is enough to show that
$\xi$ is surjective and $\eta$ is injective for
$k<i\le n+1$.

For $k<i=n+1$, $\xi$ is an isomorphism, and
for $k<i\le n$ one has $H^k(\R\Gamma_{\Phi^i}(M))=0$.

The morphism $\eta$ is injective for $i=n+1$,
and also for $k=n-1<i=n$ by \eqref{eq:k=n}.
In the remaining case $k<n-1$, 
one has $H^{k+1}(\tau^{\le n}\R\Gamma_{\Phi^{n+1}}(M))=0$.
\qed

\medskip

This $t$-structure behaves as follows under 
the local cohomology functors and the direct image functors.

\Lemma\label{lem:ld}
\bi
\item
For a locally closed subset $Z$ of $X$,
the functor $\R\Gamma_Z\cl \DC^{\rb}_\qc(\A)\to\DC^{\rb}_\qc(\A)$
is left exact, i.e.\ 
it sends $\qDA{\ge0}$ to itself.
\item
Let $j\cl U\hookrightarrow X$ be an open set of $X$.
Then $j^{-1}\cl\DC^\rb_\qc(\A)\to\DC^\rb_\qc(\A|_U)$
is exact, and
$\R j_*\cl\DC^\rb_\qc(\A|_U)\to\DC^\rb_\qc(\A)$
is left exact.
\ei
\enlemma
\proof
(i) follows immediately from
$\R\Gamma_{\Phi^n}\R\Gamma_Z(M)=\R\Gamma_Z\R\Gamma_{\Phi^n}(M)$,
which is a consequence of Lemma \ref{lem:aux} (iii), and
(ii) follows from
$\R\Gamma_{\Phi^n}\R j_*(M)=\R j_*\R\Gamma_{\Phi^n}(M)$.
\qed

\Prop\label{prop:stack}
For an open subset $U$ of $X$, let
$\CC_\qc(U)$ be the heart
$\qDA[|_U]{\le0}\cap\qDA[|_U]{\le0}$ of the $t$-structure $\qDA[|_U]{\rb}$.
Then, $U\mapsto\CC_\qc(U)$ is a stack on $X$.
\enprop
\proof
Since $X$ is Noetherian,
it is enough to show that, for open sets $U_1$, $U_2$ with $X=U_1\cup U_2$,
$M_{12}\in\CC_\qc(U_1\cap U_2)$ and $M_i\in\CC_\qc(U_i)$ 
with an isomorphism $\phi_i\cl M_i|_{U_1\cap U_2}\isoto M_{12}|_{U_1\cap U_2}$
($i=1,2$),
there exist $M\in \CC_\qc(X)$ and isomorphisms
$\psi_i\cl M|_{U_i}\isoto M_i$ ($i=1,2$)
such that the following diagram commutes.
$$\begin{CD}
M|_{U_1\cap U_2}@>{\psi_1}>> M_1|_{U_1\cap U_2}\\
@VV{\psi_2}V@VV{\phi_1}V\\
M_2|_{U_1\cap U_2}@>{\phi_2}>>M_{12}
\end{CD}
$$
Let $j_i\cl U_i\hookrightarrow X$ ($i=1,2$) and
$j_{12}\cl U_1\cap U_2\hookrightarrow X$ be the open embeddings,
and let
$$M\To[\psi_1\oplus\psi_2] \R j_{1\,*}M_1\oplus\R j_{2\,*}M_2
\To[{-\phi_1\oplus \phi_2}]
\R j_{12\,*}M_{12}\To[\ +1\ ]$$
be a distinguished triangle.
Since $\R j_{2\,*}M_2|_{U_1}\to
\R j_{12\,*}M_{12}|_{U_1}$ is an isomorphism,
$M|_{U_1}\to \R j_{1\,*}M_1|_{U_1}\simeq M_1$ is an isomorphism.
The commutativity of the diagram above is obvious by the construction.
\qed

\bigskip
Let us remark that $\bigl(\qDA{\le0}\cap \DC^{\rb}_\coh(\A),
\qDA{\ge0}\cap \DC^{\rb}_\coh(\A)\bigr)$
is not necessarily a $t$-structure on $\DC^{\rb}_\coh(\A)$
(see Example \ref{ex:Oco} and Example \ref{ex:Dco}).
In the coherent case, the proof of Theorem \ref{th:t-str} fails, because
$\R\Gamma_\Phi$ does not respect coherency.

\section{The dual standard t-structure on $\DC_\coh^\rb(\RO_X)$}
\label{sect:dual}
In the sequel, we assume that
$X$ is a finite-dimensional regular Noetherian separated scheme.
Let us denote by $\Set$ the \sd\ given by
\eq
&&\quad\ba{rl}
\Set^d\seteq&\set{Z}{\mbox{$Z$ is a closed subset of $X$ such that
$\codim Z\ge d$}}\\
=&\set{Z}{\mbox{$Z$ is a closed subset of $X$ such that
$\dim \RO_{X,x}\ge d$ for any $x\in Z$}}.\ea
\eneq

One has
\eqn
&&\sDA{\le0}\seteq\set{M\in\DC^\rb_\qc(\A)}%
{\mbox{$\Supp(H^n(M))\subset\Set^n$ for every $n\ge0$}},\\
&&\sDA{\ge0}\seteq\set{M\in\DC^{\ge0}_\qc(\A)}%
{\parbox{200pt}{$\SH^n_Z(M)=0$ for any closed subset $Z$ and $n<\codim Z$}}.
\eneqn

\Prop\label{prop:sta}
Let $Z$ be a locally closed subset of $X$.
Then $\R\Gamma_Z\cl \DC^\rb_\qc(\A)\to \DC^\rb_\qc(\A)$
is an exact functor with respect to the $t$-structure
$\sDA{\rb}$.
\enprop
\proof
It is already proved that $\R\Gamma_Z$ is left exact.
Let us prove that it is right exact,
i.e.\ it sends $\sDA{\le0}$ into itself.

We shall first prove it when $Z$ is closed.
We may assume that $X$ is affine.
Writing $Z=\cap_{i=1}^n f_i^{-1}(0)$ with $f_i\in\Gamma(X;\RO_X)$,
and using $\R\Gamma_Z=\R\Gamma_{f_1^{-1}(0)}
\circ\R\Gamma_{\cap_{i\not=1} f_i^{-1}(0)}$,
the induction on $n$ reduces to the case
$Z=f^{-1}(0)$ for some $f\in \Gamma(X;\RO_X)$.
By devissage, it is enough to prove that
$\Supp(\SH^k_Z(F[-d]))\subset \Set^{k}$ for any $k$
and $F\in\Mod_\qc(\A)$ satisfying $\Supp(F)\subset\Set^d$.
Since $\SH^k_Z$ commutes with filtrant inductive limits,
we may assume that $F$ is coherent.
Set $S\seteq\Supp(F)\in\Set^d$.
One has $\SH^k_Z(F)=0$ for $k\not=0$, $1$.
Since $\SH^0_Z(F)\subset F$,
one has $\Supp(\SH^0_Z(F))\subset \Supp(F)\in\Set^d$.
We can divide
$S=S_1\cup S_2$ with closed subsets $S_1$ and $S_2$ such that
$S_1\subset f^{-1}(0)$
and $\codim (S_2\cap f^{-1}(0))>d$.
Then on has $\Supp(\SH^1_Z(F))\subset S_2\in\Set^{d+1}$.

\smallskip
For a locally closed $Z$, let us show
$\R\Gamma_{Z}(M)\in \sDA{\le0}$ for any $M\in \sDA{\le0}$.
Writing
$Z=Z_1\setminus Z_2$ with closed $Z_2\subset Z_1\subset X$,
by the distinguished triangle
$$\R\Gamma_{Z_2}(M)\to\R\Gamma_{Z_1}(M)\to \R\Gamma_{Z}(M)\To[+1],$$
the desired result $\R\Gamma_{Z}(M)\in\sDA{\le0}$
follows from
$\R\Gamma_{Z_2}(M)[1]$, $\R\Gamma_{Z_1}(M)\in\sDA{\le0}$.
\qed

\Prop
Let $j\cl U\hookrightarrow X$ be an open embedding.
Then 
$\R j_*\cl \DC^\rb_\qc(\A|_U)\to \DC^\rb_\qc(\A)$ 
is exact with respect to
the $t$-structures $\sDA{\rb}$
and $\sDA[|_U]{\rb}$.
\enprop
\proof
Since the left exactitude has been proved, 
let us show that
$\R j_*$ sends $\sDA[|_U]{\le0}$ to $\sDA{\le0}$.
Let us denote by ${}^\Set\tau^{\le0}$ 
the truncation functor with respect to the $t$-structure $\sDA{\rb}$.
For $N\in \sDA[|_U]{\le0}$, one has
$j^{-1}({}^\Set\tau^{\le0}\R j_*N)\simeq N$, and hence
$\R j_*N\simeq \R\Gamma_U({}^\Set\tau^{\le0}\R j_*N)$
belongs to $\sDA{\le0}$ by Proposition \ref{prop:sta}.
\qed

\medskip
In the rest of this section, we treat the case where $\A=\RO_X$.
Set 
\eqn
&&\sDOc{\le0}\seteq\sDO{\le0}\cap
\DC_\coh^{\rb}(\RO_X)\quad\mbox{and}\quad
\sDOc{\ge0}\seteq\sDO{\ge0}\cap
\DC_\coh^{\rb}(\RO_X)\,,
\eneqn
and $\sDOc{\rb}\seteq(\sDOc{\le0},\sDOc{\ge0})$.

\Prop\label{prop:dualt}
\bi
\item
The pair $\sDOc{\rb}$
is a $t$-structure on $\DC_\coh^{\rb}(\RO_X)$.
\item
The equivalence of triangulated categories
$$\R\hom_{\RO_X}(-,\RO_X)
\cl \DC_\coh^{\rb}(\RO_X)\to \DC_\coh^{\rb}(\RO_X)^\op$$
sends $(\DC_\coh^{\le0}(\RO_X),\DC_\coh^{\ge0}(\RO_X))$
to $(\sDOc{\ge0}^\op,\sDOc{\le0}^\op)$.
\ei
\enprop
\begin{remark}
In  \cite[Exercise X.2]{KS} and Theorem \ref{th:tOc},
more general results are given.
\end{remark}

\proof
By Lemma \ref{lem:equi},
it is enough to show that
the functor $\R\hom_{\RO_X}(-,\RO_X)$
sends $\DC_\coh^{\le0}(\RO_X)$ and $\DC_\coh^{\ge0}(\RO_X)$
to $\sDOc{\ge0}$ and $\sDOc{\le0}$, respectively.

First assume that $M\in\DC_\coh^{\le0}(\RO_X)$.
Then one has
\eqn&&
\R\Gamma_{\Set^d}\R\hom_{\RO_X}(M,\RO_X)
\simeq \R\hom_{\RO_X}(M,\R\Gamma_{\Set^d}(\RO_X)).
\eneqn
Since $\R\Gamma_{\Set^d}(\RO_X)\in \DC_\qc^{\ge d}(\RO_X)$,
one has $\R\Gamma_{\Set^d}\R\hom_{\RO_X}(M,\RO_X)\in \DC_\qc^{\ge d}(\RO_X)$.

Next assume $M\in\DC_\coh^{\ge0}(\RO_X)$, and set
$S\seteq\Supp(H^d(\R\hom_{\RO_X}(M,\RO_X)))$.
For $x\in S$,
one has
$H^d(\R\hom_{\RO_X}(M,\RO_X))_x\simeq
\Ext^d_{\RO_{X,x}}(M_x,\RO_{X,x})\not=0$.
Since the homological dimension of $\RO_{X,x}$ is equal to
$\dim \RO_{X,x}$, one has $d\le\dim \RO_{X,x}$.
Hence one concludes that $S\in\Set^d$.
\qed

\medskip
The category $\sDO{\ge0}$ is described by a more familiar notion: 
flatness.

\Def
An object of $\DC^\rb(\RO_X)$ is called {\em with flat dimension $\le 0$}
if it is isomorphic to a bounded complex $M^\bullet$ of flat 
$\RO_X$-modules such that $M^n=0$ for $n<0$.
\edf

\Prop\label{prop:flat}
For $M \in \DC^\rb_\qc(\RO_X)$,
the following conditions are equivalent.
\bi
\item
$M\in \sDO{\ge0}$.
\item
$M$ is of flat dimension $\le0$.
\item
For any coherent $\RO_X$-module $F$, one has
$F\Ltens_{\RO_X}M\in\DO[\ge0]$.
\item
For any $N\in \DC^{\ge0}_\qc(\RO_X)$, one has
$N\Ltens_{\RO_X}M\in\DO[\ge0]$.
\ei
\enprop
\proof
The equivalence (ii)---(iv) is more or less well-known.

\smallskip
\noindent
(iv)$\Rightarrow$ (i)\quad
For every $n\ge0$, one has
$\R\Gamma_{{\Set^n}}(\RO_X)\in\DO[\ge n]$. Hence
(iv) implies that $\R\Gamma_{{\Set^n}}(M)
\simeq \R\Gamma_{{\Set^n}}(\RO_X)\Ltens_{\RO_X}M$
belongs to $\DO[\ge n]$.

\smallskip
\noindent
(i)$\Rightarrow$ (iii)\quad
For a coherent $\RO_X$-module $F$, set 
$N=\R\hom_{\RO_X}(F,\RO_X)\in \sDO{\le0}$.
Then one has
$F\Ltens_{\RO_X}M\simeq
\R\hom_{\RO_X}(N,M)\in \DO[\ge0]$.
\qed

\begin{remark}
Let $\CC_\qc\seteq\sDO{\le0}\cap\sDO{\ge0}$
be the heart of the $t$-structure $\sDO{\rb}$, and let
$\CC_\coh\seteq\CC_\qc\cap \DC_\coh^{\rb}(\RO_X)$
be the heart of the $t$-structure $\sDOc{\rb}$.
The abelian category $\CC_\coh$ is equivalent to $\Mod_\coh(\RO_X)^\op$.
It is well-known that the category $\Mod_\qc(\RO_X)$ is equivalent to
the category $\Ind(\Mod_\coh(\RO_X))$
of ind-objects of the category of coherent modules.
I conjecture that $\CC_\qc$ is equivalent to
$\Ind(\CC_\coh)\simeq \bigl(\Pro(\Mod_\coh(\RO_X))\bigr)^\op$.
Here $\Pro$ means the category of pro-objects.
\end{remark}

\section{t-structures on $\DC^\rb_\coh(\RO)$}\label{sect:tsto}

For a \sd\ $\Phi$, one sets
$$\mbox{$\qDOc{\le0}=\qDO{\le0}\cap\DOc$ and 
$\qDOc{\ge0}=\qDO{\ge0}\cap\DOc$},$$
and $\qDOc{\rb}\seteq(\qDOc{\le0},\qDOc{\ge0})$.
As seen in the example below, $\qDOc{\rb}$ is not necessarily a
$t$-structure on $\DOc$.
In this section, we give a criterian for
$\qDOc{\rb}$ to be a
$t$-structure on $\DOc$.

For an integer $n$ and $M\in\DC^\rb_\qc(\A)$, we denote 
by ${}^\Phi H^n(M)\in\qDA{\le0}\cap\qDA{\ge0}$
the $n$-th cohomology with respect to the $t$-structure
$\qDA{\rb}$.

\begin{example}\label{ex:Oco}
Let $k$ be a field, and $A=k[x]$ with an indeterminate $x$.
Let $X=\Spec(A)$ be the line.
Set
\eqn&&
\Phi^i=\left\{
\ba{ll}
\Set^0&\mbox{for $i\le0$,}\\
\Set^1&\mbox{for $i=1$, $2$,}\\
\{\emptyset\}&\mbox{for $i\ge3$.}\\
\ea\right.
\eneqn
Then, the corresponding $t$-structure on $\DO$ is given by:
\eqsub
\qDO{\le0}&=&\set{M\in \DC^{\le2}_\qc(\RO_X)}%
{\parbox{170pt}%
{$\Gamma(X;H^1(M))$ and $\Gamma(X;H^2(M))$ are torsion $A$-modules
}},\\[10pt]
\qDO{\ge0}&=&\set{M\in \DC^{\ge0}_\qc(\RO_X)}%
{\R\Gamma_{\Set^1}(M)\in \DC^{\ge2}_\qc(\RO_X)}.
\eneqsub
Let $\CC\seteq\qDO{\le0}\cap \qDO{\ge0}$
be the heart of this $t$-structure.
Since any object of $\DC^{\rb}_\qc(\RO_X)$ has a form
$\oplus M_n[-n]$ with $M_n\in\Mod_\qc(\RO_X)$,
one has
$$\CC=\set{L\oplus M[-2]}{\mbox{$L$ is a vector space over $k(x)$
and $M$ is a torsion $k[x]$-module}}.
$$
Here and in the rest of this example, we confuse $A$-modules
and quasi-coherent $\RO_X$-modules.
Since $\Hom_{\CC}(L, M[-2])=\Hom_{\CC}(M[-2],L)=0$ for
such an $L$ and $M$,
the abelian category $\CC$ is the direct sum of $\Mod(k(x))$ and
the abelian category $\Mod_\tor(k[x])$
of torsion $k[x]$-modules. 

Let $K\in\Mod(\RO_X)$ be the sheaf of rational functions.
Then $K$ and $(K/\RO_X)[-2]$ belong to $\CC$.
Hence the distinguished triangle
$$(K/\RO_X)[-2][1]\to\RO_X\to K\To[+1]$$
implies
\eqn
&&{}^\Phi H^n(\RO_X)=
\left\{
\ba{ll}
(K/\RO_X)[-2]&\mbox{for $n=-1$,}\\
K&\mbox{for $n=0$,}\\
0&\mbox{for $n\not=-1$, $0$.}\\
\ea
\right.
\eneqn
Hence, $(\qDO{\le0}\cap\DC^{\rb}_\coh(\RO_X),
\qDO{\ge0}\cap\DC^{\rb}_\coh(\RO_X))$
is {\em not} a $t$-structure on $\DC^{\rb}_\coh(\RO_X)$.

The categories
$\DC^\rb(\CC)$ and $\DC^{\rb}_\qc(\RO_X)$ are not equivalent as categories.
Indeed, the object
$K\in\CC\subset \DC^\rb(\CC)$ is a non-zero object of $\DC^\rb(\CC)$
such that
any non-zero morphism $W\to K$ in  $\DC^\rb(\CC)$
has a section.
However there is no object of $\DC^{\rb}_\qc(\RO_X)$ 
with such properties. If $M\in \DC^{\rb}_\qc(\RO_X)$ has such properties,
then, since $\Hom_{\DC^{\rb}_\qc(\RO_X)}(\RO_X[n],M)\not=0$ for some $n$,
$M$ must be a direct summand of $\RO_X[n]$. 
Since $\End_{\DC^{\rb}_\qc(\RO_X)}(\RO_X[n])\simeq k[x]$
does not have a non trivial idempotent,
$M$ must be isomorphic to $\RO_X[n]$. However,
$x\cl\RO_X[n]\to\RO_X[n]$ is a non-zero morphism but not an epimorphism.
\end{example}

In \cite{Gap}, Y. T. Siu and G. Trautmann studied
the vanishing and the coherency of local cohomologies.
Although they discussed in the analytic framework,
their main results, in our context, may be stated as follows.
Let $M\in\Mod_{\coh}(\RO_X)$,
$Z$ a closed subset of $X$ and $n$ an integer. Then one has
\bi
\item
$\SH^k_Z(M)=0$ for any $k<n$\\
$\Longleftrightarrow$
$\codim\bigl(Z\cap \Supp(\ext^k_{\RO_X}(M,\RO_X))\bigr)\ge k+n$ for any $k$.

\item
$\SH^k_Z(M)$ is coherent for any $k<n$\\
$\Longleftrightarrow$
$\codim\left(Z\cap \ol{\Supp(\ext^k_{\RO_X}(M,\RO_X))\cap(X\setminus Z)}
\right)\ge k+n$ for any $k$.
\ei

We shall generalize these statements to the derived category case.

We keep the notation: $X$ is a finite-dimensional regular 
Noetherian separated scheme.

For $M\in \DC^\rb_\coh(\RO_X)$, let us denote by
$M^*$ its dual: $M^*\seteq\R\hom_{\RO_X}(M,\RO_X)$.

\Prop\label{prop:van}
Let $M\in \DC^\rb_\coh(\RO_X)$.
\bi
\item For an integer n and a closed subset $Z$ of $X$,
$\R \Gamma_Z(M)\in \DC^{\ge n}_\qc(\RO_X)$
if and only if 
$\codim\bigl(Z\cap \Supp(H^k(M^*)\bigr)\ge k+n$ for any $k$.
\item For an integer n and a family of supports $\Phi$,
$\R \Gamma_\Phi(M)\in \DC^{\ge n}_\qc(\RO_X)$
if and only if 
$\Phi\cap \Supp(H^k(M^*))\subset\Set^{k+n}$ for every $k$.
\ei
\enprop

\proof
(ii) is a consequence of (i).
Let us show (i).
$\R \Gamma_Z(M)\in \DC^{\ge n}_\qc(\RO_X)$
if and only if 
\eq
&&\R\hom(F,M)\in \DC^{\ge n}_\coh(\RO_X)\label{eq:fm}
\eneq
for any $F\in \Mod_\coh(\RO_X)$ with $\Supp(F)\subset Z$.
Since $\bigl(\R\hom(F,M)\bigr)^*=F\Ltens_{\RO_X}M^*$,
\eqref{eq:fm} is equivalent to
$\Supp(H^k(F\Ltens_{\RO_X}M^*))\subset \Set^{k+n}$ for every $k$
by Proposition \ref{prop:dualt} (ii).
The last condition is equivalent to
$Z\cap\Supp(H^k(M^*))\subset \Set^{k+n}$
by the lemma below.
\qed

\Lemma
Let $M\in\DOc$ and $Z$, $S$ closed subsets.
Then the following conditions are equivalent.
\bi
\item
$Z\cap \Supp(\tau^{\ge0}M)\subset S$.
\item
$\Supp(\tau^{\ge0}(F\Ltens_{\RO_X}M))\subset S$
for any $F\in\Mod_\coh(\RO_X)$ with $\Supp(F)\subset Z$.
\ei
\enlemma
\proof

\noindent
(i)$\Rightarrow$(ii)\quad
For $k\ge0$, an exact sequence
$H^k(F\Ltens_{\RO_X}\tau^{<k}M)\to
H^k(F\Ltens_{\RO_X}M)\to
H^k(F\Ltens_{\RO_X}\tau^{\ge k}M)$
and $H^k(F\Ltens_{\RO_X}\tau^{<k}M)=0$
implies
$$\Supp(H^k(F\Ltens_{\RO_X}M))\subset
\Supp(H^k(F\Ltens_{\RO_X}\tau^{\ge k}M)\subset
\Supp(F)\cap\Supp(\tau^{\ge k}M)
\subset S.$$

\noindent
(ii)$\Rightarrow$(i)\quad
Let $x\in Z\cap \Supp(\tau^{\ge0}M)$. 
Take the largest $k\ge0$ such that $x\in \Supp(H^k(M))$.
If one chooses a coherent $\RO_X$-module $F$ with $\Supp(F)=Z$, then
$F_x\not=0$ and $H^k(M)_x\not=0$, which implies that
$F_x\otimes _{\RO_{X,x}}H^k(M)_x\not=0$.
On the other hand, one has
$$H^k(F\Ltens_{\RO_X}M)_x\simeq H^k(F_x\Ltens_{\RO_{X,x}}M_x)\simeq 
F_x\otimes_{\RO_{X,x}}H^k(M)_x.$$
Hence one obtains $x\in\Supp(H^k(F\Ltens_{\RO_X}M))\subset S$.
\qed

\Prop\label{prop:coh}
Let $M\in \DC^\rb_\coh(\RO_X)$,
$\Phi$ a family of supports and $n$ an integer.
Then, $\tau^{<n}(\R\Gamma_\Phi(M))\in \DC^{\rb}_\coh(\RO_X)$
if and only if 
$\Supp(H^k(M^*))\subset\Phi\cup\Psi^{k+n}$ for any $k$.
Here $\Psi^k$ is the family of supports consisting
of closed subsets $Z$ such that $\Phi\cap Z\subset\Set^k$.
\enprop
\proof
Assume first $M'\seteq
\tau^{<n}(\R\Gamma_\Phi(M))\in \DC^{\rb}_\coh(\RO_X)$.
Let us complete the morphism $M'\to M$ to a distinguished triangle
$M'\to M\to M''\To[+1]$.
Since one has
$\R\Gamma_\Phi(M'')\in \DC^{\ge n}_\qc(\RO_X)$,
Proposition \ref{prop:van}
implies that $\Supp(H^k(M''^*))\subset\Psi^{k+n}$ for every $k$.
The distinguished triangle
$M''^*\to M^*\to M'^*\To[+1]$ implies that
$\Supp(H^k(M^*))\subset\Supp(H^k(M'^*))\cup\Supp(H^k(M''^*))$.
Since  $\Supp(H^k(M'^*))\subset\Phi$ and $\Supp(H^k(M''^*))\subset\Psi^{k+n}$,
we obtain $\Supp(H^k(M^*))\subset\Phi\cup\Psi^{k+n}$.

\smallskip
Conversely, assuming that
$\Supp(H^k(M^*))\subset\Phi\cup\Psi^{k+n}$ for every $k$,
we shall prove $\tau^{<n}(\R\Gamma_\Phi(M))\in\DOc$.
By devissage, one may assume that
$M^*=F[-k]$ for $F\in\Mod_\coh(\RO_X)$ with
$\Supp(F)\subset\Phi\cup\Psi^{k+n}$.
Then there is an exact sequence
$0\to F'\to F\to F''\to 0$
with
$F'$, $F''\in\Mod_\coh(\RO_X)$
sucth that
$\Supp(F')\subset\Psi^{k+n}$ and $\Supp(F'')\subset\Phi$.
Set $M'\seteq (F'[-k])^*$ and $M''\seteq (F''[-k])^*$.
Then $\R\Gamma_\Phi(M')\in\DO[\ge n]$ by Proposition \ref{prop:van}.
Since $\R\Gamma_\Phi(M'')\simeq M''$, one has a distingushed triangle
$M''\to \R\Gamma_\Phi(M)\to \R\Gamma_\Phi(M')\To[+1]$.
Since $\tau^{<n}\R\Gamma_\Phi(M')=0$, one has
$\tau^{<n}\R\Gamma_\Phi(M)\simeq\tau^{<n}M''\in \DC^\rb_\coh(\RO_X)$.
\qed

\medskip
It is sometimes convenient to use 
$\Z$-valued functions on $X$ instead of \sds.
We say that a bounded
$\Z$-valued function $p$ on $X$ is a {\em supporting function}
if $p(y)\ge p(x)$ whenever $y\in\ol{\{x\}}$.

\Lemma
By the following correspondence, the set of \sds\ 
is isomorphic to
the set of supporting functions.
To a \sd\ $\Phi$, one associates
the supporting function $p_\Phi$
given by
$p_\Phi(x)\seteq\max\set{n\in\Z}{\ol{\{x\}}\in\Phi^n}$.

Conversely, to a supporting function $p$ on $X$,
one associates the \sd\ $\Phi_p$ given by
$\Phi_p^n=\set{Z}{\mbox{$Z$ is a closed subset such that
$p(z)\ge n$ for any $z\in Z$}}$.
\enlemma
This lemma immediately follows from 
the fact that any closed subset is a union of finitely many
irreducible closed subsets and that any irreducible subset 
has a generic point.

\medskip

One has
\eqn
&&\ba{rl}
p_{\tv}(x)&=0,\\
p_\Set(x)&=\codim(\ol{\{x\}})=\dim(\RO_{X,x}),\\[2pt]
p_{\sigma^{\le n}\Phi}(x)&=\min(p_\Phi(x),n),
\ea\qquad\ba{ll}
p_{\Phi\cap\Psi}(x)&=\min(p_{\Phi}(x),p_{\Psi}(x)),\\[2pt]
p_{\Phi\cup\Psi}(x)&=\max(p_{\Phi}(x),p_{\Psi}(x)).\\[2pt]
  & 
\ea
\eneqn

\medskip
For two \sds\ $\Phi$ and $\Psi$, we define the
\sd\ $\Phi\circ\Psi$ by
\eq
(\Phi\circ\Psi)^n=\bigcup_{n=i+j}(\Phi^i\cap\Psi^j)
\eneq

Note that
$\circ$ is commutative and associative, and
$\tv$ is the unit with respect to $\circ$:
$\tv\circ\Phi=\Phi$ for every $\Phi$.

The following lemma is obvious.

\Lemma\label{lem:circ1}
Let $\Phi$ and $\Psi$ be \sds.
\bi
\item
one has
$$\bigcup_{n=i+j}(\Phi^i\cap\Psi^j)=
\bigcap_{n+1=i+j}(\Phi^i\cup\Psi^j).$$
\item
Let $Z$ be an irreducible closed subset of $X$ such that
$Z\in\Phi^a$ and $Z\not\in\Phi^{a+1}$.
Then $Z\in\Psi^\rb$ if and only if $Z\in (\Phi\circ\Psi)^{a+b}$.
\item
$p_{\Phi\circ\Psi}(x)=p_\Phi(x)+p_{\Psi}(x)$ for any $x\in X$.
\ei
\enlemma
\proof
(iii) is obvious. Let us show (i) and (ii).
The inclusion $\bigcup_{n=i+j}(\Phi^i\cap\Psi^j)\subset
\bigcap_{n+1=i+j}(\Phi^i\cup\Psi^j)$
is obvious.
Hence it is enough to show that
for any  irreducible closed subset $Z$, the conditions
$Z\in\Phi^a$, $Z\not\in\Phi^{a+1}$ and
$Z\in \bigcap_{a+b+1=i+j}(\Phi^i\cup\Psi^j)$ imply
$Z\in\Psi^\rb$.
Since $Z\in\Phi^{a+1}\cup\Psi^\rb$, one obtains
$Z\in\Psi^\rb$.
\qed

\Lemma\label{lem:circ2}
Let $\Phi$, $\Psi_1$ and $\Psi_2$ be three \sds.
\bi
\item
$\Phi\circ(\Psi_1\cap\Psi_2)=(\Phi\circ\Psi_1)\cap(\Phi\circ\Psi_2)$
and
$\Phi\circ(\Psi_1\cup\Psi_2)=(\Phi\circ\Psi_1)\cup(\Phi\circ\Psi_2)$.
\item
If $\Phi\circ\Psi_1\subset\Phi\circ\Psi_2$,
then one has $\Psi_1\subset\Psi_2$.
\item
If $\Phi\circ\Psi_1=\Phi\circ\Psi_2$,
then one has $\Psi_1=\Psi_2$.
\ei
\enlemma
\proof
(i) is obvious.
(ii) follows from Lemma \ref{lem:circ1} (iii),
since $\Psi_1\subset\Psi_2$ if and only if
$p_{\Psi_1}(x)\le p_{\Psi_2}(x)$
for every $x\in X$. 
(iii) is an immediate consequence of (ii).
\qed

\medskip
\Lemma
Let $\Phi$ and $\Theta$ be a pair of \sds.
For an integer $a$, set
\eqn
&&
\ba{rl}
\Psi^n&\seteq\set{Z}{\mbox{$Z$ is a closed subset such that
$\Phi^k\cap Z\subset\Theta^{k+n}$ for every $k$}},\\
(\Psi_a)^n&\seteq
\set{Z}{\mbox{$Z$ is a closed subset such that
$\Phi^k\cap Z\subset\Theta^{k+n}$ for every $k\le a$}}\\
&\phantom{:}=\set{Z\in(\Psi_{a-1})^n}{\Phi^a\cap Z\subset\Theta^{a+n}}.
\ea
\label{eq:psia}
\eneqn
and $\Psi\seteq\{\Psi^n\}_{n\in\Z}$,
$\Psi_a\seteq\{\Psi_a^n\}_{n\in\Z}$.
\bi
\item
$\Psi$ and $\Psi_a$ are \sds, and
$p_\Psi(x)=\min\set{p_\Theta(y)-p_\Phi(y)}{y\in\ol{\{x\}}}$.
\item
If a \sd\ $\Psi'$ satisfies $\Phi\circ\Psi'=\Theta$, then
$\Psi'=\Psi$.
\item
There exists a \sd\ $\Psi'$ such that $\Phi\circ\Psi'=\Theta$
if and only if
$0\le p_\Phi(y)-p_\Phi(x)\le p_\Theta(y)-p_\Theta(x)$
whenever
$y\in\ol{\{x\}}$.
\item
If $\Phi\circ\Psi=\Theta$, then
one has
$(\sigma^{\le a}\Phi)\circ \Psi_a=\Theta$
for every integer $a$.
\item
Assume that $\Phi\circ\Psi=\Theta$.
Then one has
$(\Psi_{a-1})^n\subset\Phi^a\cup\Psi^n$ for every $n$.
\item
Assume that
$(\Psi_{a-1})^n\subset\Phi^a\cup(\Psi_a)^n$
for every $a$ and $n$.
Then $\Phi\circ \Psi=\Theta$.
\ei
\enlemma
\proof
(i) is obvious.

\smallskip
\noindent
(ii)\quad
One has
$\Psi'\subset\Psi$ and $\Phi\circ\Psi\subset\Theta$.
Hence
$$\Theta=\Phi\circ\Psi'\subset\Phi\circ\Psi\subset\Theta.$$
Hence $\Phi\circ\Psi'=\Phi\circ\Psi$
and Lemma \ref{lem:circ2} (iii) implies $\Psi'=\Psi$.

\smallskip
\noindent
(iii)  is obvious.

\smallskip
\noindent
(iv)\quad
We shall apply (iii) and (ii).
If $y\in\ol{\{x\}}$, then one has
$$p_{\sigma^{\le a}\Phi}(y)-p_{\sigma^{\le a}\Phi}(x)
=\min(p_\Phi(y),a)-\min(p_\Phi(x),a)
\le p_\Phi(y)-p_\Phi(x)\le p_\Theta(y)-p_\Theta(x).$$
Here the first inequality follows from $p_\Phi(y)\ge p_\Phi(x)$.

\smallskip
\noindent
(v)\quad
Let $Z\in(\Psi_{a-1})^n$ be an irreducible closed subset.
We shall show that if $Z\not\in\Phi^a$,
then $Z\in\Psi^n$.
Let us take an integer $i$ such that
$Z\in\Phi^i$ and $Z\not\in\Phi^{i+1}$.
Then one has $i<a$.
Since $Z\in(\Psi_{a-1})^n$, one has
$Z=\Phi^i\cap Z\subset\Theta^{i+n}=(\Phi\circ\Psi)^{i+n}$.
Hence $Z\in\Psi^n$.

\smallskip
\noindent
(vi)\quad
Let $Z\in\Theta^n$ be an irreducible 
closed subset.
Let us show
$Z\in(\Phi\circ\Psi)^n$.
Take an integer $i$
such that $Z\in\Phi^i$ and $Z\not\in\Phi^{i+1}$.
Let us show $Z\in\Psi^{n-i}$.
Since $\Psi_a=\Psi$ for $a\gg0$,
it is enough to show
$Z\in (\Psi_a)^{n-i}$ for every $a$.
We shall show it by the induction on $a$.
It is obvious that $Z\in (\Psi_a)^{n-i}$ for $a\ll0$.
Assuming that $Z\in(\Psi_{a-1})^{n-i}$,
let us show $Z\in(\Psi_a)^{n-i}$.
Since $Z\in (\Psi_{a-1})^{n-i}\subset\Phi^a\cup(\Psi_a)^{n-i}$,
we may assume that
$Z\in\Phi^a$, which implies $a\le i$.
Therefore one has $\Phi^a\cap Z=Z\in\Theta^n\subset\Theta^{n-i+a}$.
Together with $Z\in(\Psi_{a-1})^{n-i}$,
one obtains $Z\in (\Psi_a)^{n-i}$.
\qed

\medskip
For a \sd\ $\Phi$, let $\Phi_*$ denote the \sd\ given by
\eq
&&\Phi_*{}^n\seteq\set{Z}{\mbox{$Z$ is a closed subset such that
$Z\cap\Phi^k\subset\Set^{n+k}$ for every $k$}}.
\eneq

Now we are ready to give a criterian for
$\qDOc{\rb}$ to be a $t$-structure,
which is a generalization of
\cite[Exercise X.2]{KS}.

\Theorem\label{th:tOc}
Let $\Phi$ be a \sd.
Then the following conditions are equivalent.
\bi
\item
$\qDOc{\rb}$ is a $t$-structure.
\item
For any irreducible closed subsets
$Z$ and $S$ such that $S\subset Z$ and $S\in\Phi^k$, one has
$Z\in\Phi^{k+\codim(Z)-\codim(S)}$.
In terms of supporting functions, one has
$$\mbox{$p_\Phi(y)-p_\Phi(x)\le \codim(\ol{\{y\}})-\codim(\ol{\{x\}})
=\dim(\RO_{X,y})-\dim(\RO_{X,x})$
if $y\in\ol{\{x\}}$.}$$
\item
$\Phi\circ\Phi_*=\Set$.
\item
There exists a \sd\ $\Psi$ such that
$\Phi\circ\Psi=\Set$.
\item
$(\sigma^{<n}\Phi)_*^k\subset\Phi^n\cup
(\sigma^{\le n}\Phi)_*{}^k$
for every $n$ and every $k$.
\ei
Moreover if these conditions are satisfied,
the equivalence
$$\R\hom{\RO_X}(-,\RO_X)\cl \DOc\isoto \DOc^\op$$
sends the $t$-structure $\qDOc[\Phi]{\rb}$
to
$\qDOc[\Phi_*]{\rb}^\op$.
\entheorem

\proof
The last statement easily follows from Proposition \ref{prop:van}.
The equivalence of
(ii)---(v) is already shown.

\smallskip
Let us show (v)$\Rightarrow$ (i).
The proof is similar to
the one of Theorem \ref{th:t-str}.
It is enough to show that any $M\in\DOc$ satisfies the following property:
\eq
&&\parbox{350pt}{There exists a distinguished triangle
$M'\to M\to M''\To[+1]\,$ with
$M'\in\qDOc{<0}$ and $M''\in\qDOc{\ge0}$.}\label{eq:disto}
\eneq
For $n\in \Z$, let us consider the following statement:
\begin{equation}
\parbox{400pt}{
The property \eqref{eq:disto} holds
if $M\in\qDOc[\sigma^{\le n}\Phi]{\ge0}$.}
\tag*{(\ref{eq:disto})${}_n$}
\end{equation}
We shall prove \eqref{eq:disto}${}_{n}$ 
by assuming \eqref{eq:disto}${}_{n+1}$.
Let us consider a distinguished triangle
$$\tau^{\le n}\R\Gamma_{\Phi^{n+1}}(M)\to M\to M''\To[+1].$$
In the course of the proof of Theorem \ref{th:t-str},
we have proved that
$\tau^{\le n}\R\Gamma_{\Phi^{n+1}}(M)=\SH^n_{\Phi^{n+1}}(M)[-n]
\in\qDO{<0}$ and
$M''\in\qDOX[{\sigma^{\le n+1}\Phi}]{\ge0}$.
Hence it is enough to show
that
$\SH^n_{\Phi^{n+1}}(M)$
is coherent.
Since
$\SH^i_{\Phi^k}(M)=0$ for $i<k\le n$,
Proposition \ref{prop:van}
implies that
$\Supp(H^i(M^*))\cap\Phi^k\subset\Set^{k+i}$ for every $k\le n$ and $i$.
Hence one has
$\Supp(H^i(M^*))\in(\sigma^{\le n}\Phi)_*{}^i$ for every $i$.
Since
$(\sigma^{\le n}\Phi)_*{}^i\subset\Phi^{n+1}\cup
(\sigma^{\le n+1}\Phi)_*{}^i$,
Proposition \ref{prop:coh}
implies that
$\SH^i_{\Phi^{n+1}}(M)$ is coherent for every $i<n+1$.

\smallskip
Let us show that
(i) implies (ii).
In order to see this, it is enough to show that the following
situation cannot happen:
\eqn
&&\parbox{420pt}
{$Z$ is an irreducible closed subset of $X$ and
$S$ is an irreducible closed subset of $Z$ with codimension $1$
such that
$Z\in\Phi^a$, $Z\not\in\Phi^{a+1}$ and
$S\in\Phi^\rb$ with $b>a+1$.}
\eneqn

Set $M=\RO_Z[-a]\in\qDOc{\le0}$.
Here $\RO_Z$ is the structure sheaf of $Z$ endowed
with the reduced scheme structure.
Let
$$M'\to M\to M''\To[+1]$$
be a distinguished triangle with
$M'\in\qDOc{<0}$ and $M''\in\qDOc{\ge0}$.
Then one has $\Supp(M'')\subset Z$, and hence
$M''=\R\Gamma_{\Phi^a}(M'')\in \DOc[\ge a]$,
which implies $M'\in  \DOc[\ge a]$.
The exact sequence
$0\to H^a(M')\to\RO_Z$,
along with $\Supp(H^a(M'))\in\Phi^{a+1}$,
implies that $H^{a}(M')=0$.
Hence one has $M'\in \DOc[>a]$.
On the other hand, one has
the exact sequence
$$\SH^{a+1}_S(M')\to \SH^{a+1}_S(M)\to \SH^{a+1}_S(M'').$$
Since $\R\Gamma_S(M'')\in\DO[\ge b]$,
one has $\SH^{a+1}_S(M
'')=0$.
Since
$M'\in  \DOc[>a]$, one has
$\SH^{a+1}_S(M')=\Gamma_S(H^{a+1}(M'))$,
which is a coherent $\RO_X$-module.
Hence $\SH^{a+1}_S(M)=\SH^1_S(\RO_Z)$ is a coherent $\RO_X$-module.
This is a contradiction.
The last step follows from either Proposition \ref{prop:coh} or
the following easy lemma, whose proof is
omitted.
\qed

\Lemma
Let $A$ be a $1$-dimensional
Noetherian local ring with a maximal ideal $\m$.
Then $H^1_{\m}(A)$ is not a finitely generated $A$-module.
\enlemma

\section{t-structure on $\DC_\qc^\rb(\RD_X)$ and $\DC_\hol^\rb(\RD_X)$}%
\label{sec:D}

In the sequel, we shall treat the case $\A=\RD$.
Let $X$ be an algebraic manifold over a field $k$ of characteristic $0$, 
i.e.\ a quasi-compact separated scheme 
smooth over $k$.
Let $\RD_X$ be the sheaf of differential operators on $X$.
Let us denote by
$\DC_\qc^\rb(\RD_X)$ the derived category of
bounded complexes of $\RD_X$-modules with quasi-coherent cohomologies.
Let us denote by
$\DC_\coh^\rb(\RD_X)$, $\DC_\h^\rb(\RD_X)$, $\DC_\rh^\rb(\RD_X)$
the full subcategories of
$\DC_\qc^\rb(\RD_X)$ consisting of
bounded complexes with
coherent, holonomic, regular holonomic $\RD_X$-modules as cohomologies,
respectively.
For a morphism $f\cl X\to Y$ of algebraic manifolds,
we denote by
$\D f^*\cl\DC_\qc^\rb(\RD_Y)\to\DC_\qc^\rb(\RD_X)$ and $\D f_*\cl
\DC_\qc^\rb(\RD_X)\to\DC_\qc^\rb(\RD_Y)$ the inverse image 
and the direct image functors (see e.g.\ \cite{K1}).
Note that they respect $\DC_\h^\rb$ and $\DC_\rh^\rb$.

Let us define a $t$-structure on $\DC_\qc^\rb(\RD_X)$ as follows.
\eq
\ba{ll}
\sDD[\le0]&=\set{M\in \DC_\qc^\rb(\RD_X)}%
{\mbox{$\Supp(H^n(M))\subset\Set^n$ for every $n$}},\\[3pt]
\sDD[\ge0]&=\set{M\in \DC_\qc^\rb(\RD_X)}%
{\mbox{$\R\Gamma_{\Set^n}(M)\in\DC_\qc^{\ge n}(\RD_X)$ for every $n$}}.
\ea
\eneq
It is indeed a $t$-structure by Theorem \ref{th:t-str}.
Let $\CC_\qc=\qDle\cap \qDge$ be its heart.

We note that 
$\bigl(\qDle\cap\DC_\coh^\rb(\RD_X),\qDge\cap\DC_\coh^\rb(\RD_X)\bigr)$
is {\em not} a $t$-structure on $\DC_\coh^\rb(\RD_X)$ (when $\dim X>1$),
as seen in the following
example.

\begin{example}\label{ex:Dco}
Set $X=\Spec(\C[x,y])$ and $Y=\Spec(\C[t])$
with indeterminates $x$, $y$ and $t$.
Set $t=f(x,y)\seteq xy\cl X\to Y$.
Set $\CM\seteq\RD_X/\RD_X(x\partial_x-y\partial_y)$
and $\CN\seteq\RD_f$. Note that $\CN\simeq \RO_X\otimes_{\RO_Y}\RD_Y
=\oplus_{n=0}^\infty\RO_Xv_n$ with $v_n=\one_{X\to Y}\otimes\partial_t^n$
as an $\RO_X$-module.
The defining relations of $\{v_n\}_{n\in\Z_{\ge0}}$ as a $\RD_X$-module are
$$\partial_x v_n=yv_{n+1}\mbox{ and }\partial_y v_n=xv_{n+1}.$$
Hence $\CN\in\CC_\qc$.
One has a morphism $\CM\to \CN$ by 
$1\,\mathrm{mod}\, \RD_X(x\partial_x-y\partial_y)\mapsto v_0$.
Set $E=\C^{\oplus\Z_{>0}}=\smash{\mathop\oplus\limits_{n\in\Z_{>0}}}\C w_n$.
Then one has $\B_{\{0\}|X}\otimes_{\C}E[-2]\in\CC_\qc$.
There is an exact sequence
$$0\to \CM\to \CN\To[g]\B_{\{0\}|X}\otimes_{\C}E\to 0$$
in $\Mod(\RD_X)$.
Here $g$ is given by 
\eqn&&
g(v_n)
=\sum_{1\le k\le n}(-1)^{n-k}\dfrac{1}{(n-k)!}
(\partial_x\partial_y)^{n-k}\delta\otimes w_k\,,
\eneqn
where $\delta\in \B_{\{0\}|X}$ is the generator with the defining relations
$x\delta=y\delta=0$.
Hence one has a distinguished triangle
$\B_{\{0\}|X}\otimes_{\C}E[-2][1]\to \CM\to\CN\To[+1]\,$.
Thus we obtain
\eqn
&&{}^\Set H^n(\CM)=\left\{
\ba{ll}
\B_{\{0\}|X}\otimes_\C E[-2]&\mbox{for $n=-1$,}\\[3pt]
\CN&\mbox{for $n=0$,}\\[3pt]
0&\mbox{for $n\not=-1$, $0$.}\\
\ea\right.
\eneqn
Since $\B_{\{0\}|X}\otimes_\C E$ is not coherent, 
$\bigl(\qDle\cap\DC_\coh^\rb(\RD_X),\qDge\cap\DC_\coh^\rb(\RD_X)\bigr)$
is {\em not} a $t$-structure on $\DC_\coh^\rb(\RD_X)$.
Note that one has $\SH^0_{\Set^1}(\CM)=0$
and $\SH^1_{\Set^2}(\CM)\simeq \B_{\{0\}|X}\otimes_{\C}E$.
\end{example}

\medskip
Let us denote by
$\qDhle\seteq\qDle\cap\DC_\hol^\rb(\RD_X)$
and $\qDhge\seteq\qDge\cap\DC_\hol^\rb(\RD_X)$.
Similarly we define
$\qDrhle$ and $\qDrhge$.

\Theorem
$(\qDhle,\qDhge)$ 
and $(\qDrhle,\qDrhge)$ are
a $t$-structure on $\Dh$ and $\DC_\rh^\rb(\RD_X)$,
respectively.
\entheorem

\proof
In order to show that $(\qDhle,\qDhge)$
is a $t$-structure on $\Dh$,
it is enough to show that,
for any any $\CM\in\Dh$, there exists a distinguished triangle
$\CM'\to \CM\to \CM''\To[+1]$
with $\CM'\in\qDh{<0}$ and $\CM''\in\qDh{\ge0}$.

Let us show this by the induction on
the codimension $d$ of $S\seteq\Supp(\CM)$. 
One has $\tau^{<d}\CM\in\qDh{<0}$.
By the distinguished triangle
$\tau^{<d}\CM\to \CM\to \tau^{\ge d}\CM\To[+1]$,
we may assume that
$\CM\simeq\tau^{\ge d}\CM$ by Lemma \ref{lem:t-str}.

Let $S_0$ be a $d$-codimensional smooth open subset of $S$
such that $S_1\seteq S\setminus S_0$ is of codimension $>d$.
Set $U\seteq X\setminus S_1$. Then one has $S_0=U\cap S$.
Let $j\cl U\hookrightarrow X$ and
$i\cl S_0\hookrightarrow U$ be the open embedding and
the closed embedding, respectively.
Then there exists $\CN\in\DC_\h^{\ge d}(\RD_{S_0})$
such that $\CM |_U$ is isomorphic to $\D i_*\CN$.
By shrinking $S_0$ if necessary, we may assume that
the cohomologies of $\CN$ are locally free $\RO_{S_0}$-modules.
Hence $\D i_*\CN$ belongs to $\qDhge[U]$.
Hence $\CM''\seteq\R j_*\D i_*\CN$ belongs to $\qDhge$
by Lemma \ref{lem:ld}.
Let us consider a distinguished triangle
$$\CM'\to \CM\to \CM''\To[+1].$$ 
Since $\Supp (\CM')\subset S_1$,
the codimension of $\Supp(\CM')$ is greater than $d$.
Then the induction proceeds by Lemma \ref{lem:t-strd}.

The regular holonomic case is proved similarly,
because the regular holonomicity is also
preserved by the direct image functors.
\qed

\begin{remark}
By Proposition \ref{prop:stack},
$U\mapsto{}^\Set\DC^{\le0}_*(\RD_{U})\cap{}^\Set\DC^{\ge0}_*(\RD_{U})$
is a stack on $X$ for $*=\qc$, $\hol$, $\rh$.
\end{remark}

For a closed point $x\in X$, let us denote by
$\B_{x|X}$ the regular holonomic $\RD_X$-module $\SH^{\dim X}_{\{x\}}(\RO_X)$.

\Prop\label{lem:flat}
For $\CM\in\Dh$, the following conditions are equivalent.
\bi
\item $\CM\in\qDhge$.
\item $\CM$ is of flat dimension $\le0$ as a complex of $\RO_X$-modules.
\item For any closed subset $Z$ of $X$, one has $\SH^i_Z(\CM)=0$ 
for any $i<\codim Z$.
\item For any locally closed subset $Z$ of $X$, one has 
$\SH^i_Z(\CM)=0$ for any $i<\codim Z$.
\item $\Ext^i_{\RD_X}(\B_{x|X},\CM)=0$ for any closed point
$x\in X$ and $i<\dim X$.
\item $H^i(\ko_x\Ltens_{\RO_X}\CM)=0$ for any closed point 
$x\in X$ and $i<0$.
Here $\ko_x$ is the sheaf $\RO_{X}/\m_x$ with the ideal $\m_x$ 
of functions vanishing at $x$.
\ei
\enprop
\proof
The equivalence of (i)---(iv) are evident.
Since $\R\hom_{\RD_X}(\B_{x|X},\CM)\simeq
\ko_x\Ltens_{\RO_X}\CM[-\dim X]$, the implication
(i)$\Rightarrow$ (v) and the equivalence (v)$\Leftrightarrow$(vi)
are evident.

It remains to prove
(vi)$\Rightarrow$ (i).
Let us assume that $\CM$ satisfies (vi).
Let us show $\CM\in\qDhge$ by the induction on
the codimension $d$ of $S\seteq\Supp(\CM)$.
Let $S_0$ be a $d$-codimensional smooth open subset of $S$
such that $S_1\seteq S\setminus S_0$ is of codimension $>d$.
Let $i\cl S_0\to X$ be the inclusion.
Then there exists $\CN\in\Dh[S_0]$
such that $\CM$ is isomorphic to $\D i_*\CN$ on a neighborhood of $S_0$.
By shrinking $S_0$ if necessary, we may assume that
the cohomologies of $\CN$ are locally free $\RO_{S_0}$-modules.
Then $\ko_x\otimes_{\RO_{S_0}}\CN\simeq \ko_x\otimes_{\RO_X}\CM[d]$,
and hence one has $H^k(\CN)_x=0$ for any $k<d$ and any
closed point $x$ of $S_0$.
Therefore, $\CN\in\DC^{\ge d}(\RD_{S_0})$
and $\CM''\seteq\D i_*\CN\simeq\R \Gamma_{S_0}\CM$ belongs to
$\qDhge$.
Consider a distinguished triangle
$$\CM'\to\CM\to \CM''\To[+1].$$
The exact sequence $H^{k-1}(\ko_x\Ltens_{\RO_X}\CM'')
\to H^{k}(\ko_x\Ltens_{\RO_X}\CM')\to H^{k}(\ko_x\Ltens_{\RO_X}\CM)$
implies that $H^{k}(\ko_x\Ltens_{\RO_X}\CM')=0$ for $k<0$.
Since $\codim \Supp(\CM')>d$, the induction hypothesis
implies that $\CM'\in\qDhge$. Hence one concludes that
$\CM\in\qDhge$.
\qed

\medskip
Let us assume that the base field $k$ is the complex number field $\C$,
and let us denote by $X_\an$ the associated complex manifold.
Let us denote by $\DC_c^\rb(\C_{X_\an})$ the derived category
of bounded complexes of $\C_{X_\an}$-modules with constructible cohomologies.

\Theorem
The equivalence of triangulated categories
$$\Sol_X\seteq\R\hom_{\RD_X}(-, \RO_{X_\an})\cl\Drh\isoto
\DC_c^\rb(\C_{X_\an})^\op$$ sends the $t$-structure
$(\qDrhle,\qDrhge)$ to
$(\DC_c^{\ge0}(\C_{X_\an})^\op,\DC_c^{\le0}(\C_{X_\an})^\op)$.
\entheorem
\proof
Since $\R\hom_{\RD_X}(-, \RO_{X_\an})$ is an equivalence 
of triangulated categories,
it is enough to show
that $\CM\in \Drh$ belongs to
$\qDrhge$ if and only if $\R\hom_{\RD_X}(\CM, \RO_{X_\an})$
belongs to $\DC_c^{\le0}(\C_{X_\an})$.
This immediately follows from Proposition \ref{lem:flat} and
$$\R\hom_{\RD_X}(\CM,\RO_{X_\an})_x\simeq\Hom_\C(\C_x\otimes_{\RO_X}\CM,\C).$$
\qed

\section{Proof of Lemma \ref{lem:indinj}}\label{sec:proof}

In this section, we give a proof of Lemma \ref{lem:indinj}.

\Lemma\label{lem:arb}
Let $\M$ be a coherent $\A$-module,
and $\N$ a {\rm(}not necessarily coherent\/{\rm)} $\A$-submodule of $\M$,
and $Z$ a closed subset of $X$. Then, for the
generic point $\xi$ of any irreducible component of $Z$,
there exist a coherent $\A$-submodule $\L$ of $\M$ and an open neighborhood
$U$ of $\xi$ such that
$$\N|_U=(\N_{X\setminus Z}+\L)|_U.$$
\enlemma
\proof
Assume first that $X$ is affine.
Set $R=\RO_X(X)$, $A=\A(X)$, $M=\M(X)$.
Let $L\subset M$ be the inverse image of
$\N_\xi$ by the morphism
$M\to\M_\xi$. Then $L$ is a finitely generated $A$-module.
Let $\L$ be the coherent $\A$-module associated with $L$:
$\L\seteq\A\otimes_AL\simeq\RO_X\otimes_RL$.
Then one has
$\N_\xi=\L_\xi$.
Hence there exists an open neighborhood $U$ of $\xi$ such that
$\L|_U\subset \N|_U$.
We may assume further that $U\cap Z=U\cap \ol{\{\xi\}}$.
Let us show that
$\N|_U=(\N_{X\setminus Z}+\L)|_U$. It is evident that
the last member is contained in the first.
We shall prove the opposite inclusion.
Let $V$ be an affine open subset contained in $U$.
Then $\M(V)=\RO_X(V)\otimes_R M$.
Let $N_V$ be the inverse image of $\N(V)$
by the morphism $M\to\M(V)$.
Then one has
$\N(V)=\RO_X(V)\otimes_RN_V$.

Assume first $\xi\in V$.
Then by the commutative diagram
$$\begin{CD}
N_V@>>>\N(V)@>>>\N_\xi\\
@VVV@VVV@VVV\\
M@>>>\M(V)@>>>\M_\xi
\end{CD}$$
$N_V$ is contained in $L$.
Hence $\N(V)\simeq\RO_X(V)\otimes_RN_V
\subset \RO_X(V)\otimes_RL\simeq \L(V)$.
If $\xi$ is not contained in $V$, then one has
$V\subset X\setminus Z$, and hence
$\N(V)\subset (\N_{X\setminus Z})(V)$.
Hence one obtains $\N|_U\subset(\N_{X\setminus Z}+\L)|_U$.

In the general case, let us take an affine open subset $W$ of $X$
containing $\xi$. Then it is enough to remark that any 
coherent $\A|_W$-submodule
of $\M|_W$ can be extended to a coherent $\A$-submodule of $\M$.
Note that a quasi-coherent $\A$-submodule of a coherent $\A$-module 
is coherent.
\qed

The next lemma is not used in this paper.
\Lemma
Let $\M$ be a coherent $\A$-module,
and $\N$ a {\rm(}not necessarily coherent\/{\rm)} $\A$-submodule of $\M$.
Then, there exist finite families of open subsets $U_i$ of $X$ and 
coherent $\A$-submodules $\M_i$ of $\M$ 
such that
$$\N=\sum_i(\M_i)_{U_i}.$$
\enlemma
\proof
Let $\mathcal{W}$ be the set of open subsets $W$ of
$X$ such that there exist finite families of open subsets $U_i$ of 
$W$ and
coherent $\A$-submodules $\M_i$ of $\M$ such that
$\N|_W=(\sum_i(\M_i)_{U_i})|_W$.
Since $X$ is Noetherian, $\mathcal{W}$ has a maximal element.
Let $W$ be such a maximal element.
Assuming that $W\not=X$, let us derive a contradiction.
Let $\xi$ be the generic point of an irreducible component of
$X\setminus W$.
Then by Lemma \ref{lem:arb}, there exist an open neighborhood $U$
of $\xi$ and a coherent $\A$-submodule $\L$ such that
$\N|_U=(\N_W+\L)|_U$.
Hence one has 
$\N|_{W\cup U}=\Bigl(\sum_i(\M_i)_{U_i}+\L_U\Bigr)|_{W\cup U}$.
This contradicts the choice of $W$.
\qed

\medskip
The following lemma,
an analogue of \cite{Gol},
is a corollary of Lemma \ref{lem:arb},

\Lemma
Let $\CM$ be an $\A$-module.
Then the following conditions are equivalent.
\bi
\item
$\CM$ is an injective $\A$-module.
\item
For any coherent $\A$-module $\CF$,
one has $\ext_\A^1(\CF,\CM)=0$ and
$\hom_{\A}(\CF,\CM)$ is a flabby sheaf.
\ei
\enlemma
\proof
(i)$\Rightarrow$(ii) is well-known.
Let us show (ii)$\Rightarrow$(i).
In order to see that $\M$ is injective, it is enough to show that
for any left ideal $\I$ of $\A$ with $\J\not=\A$
and a morphism $\phi\cl \I\to \M$, there exists
a left ideal $\I'$ strictly containing $\I$ such that 
$\phi$ extends to a morphism $\I'\to\M$.
Let $Z=\Supp(\A/\I)\not=\emptyset$, and let $\xi$ be the generic point 
of an irreducible component of $Z$.
Then by Lemma \ref{lem:arb},
there exist a neighborhood $U$ of $\xi$ and
a coherent left ideal $\J$ of $\A$ such that
$\I|_U=(\A_{X\setminus Z}+\J)|_U$.
Set $\I'=\I+\A_U\not=\I$.
The exact sequence
$0\to \I\to\I' \to (\A/\J)_{Z\cap U}\to 0$
induces an exact sequence
$$\Hom_\A(\I',\M)\to \Hom_\A(\I,\M)\to \Ext^1_\A((\A/\J)_{Z\cap U},\M).$$
Here the last term vanishes by the assumptions:
\eqn
\Ext^1_\A((\A/\J)_{Z\cap U},\M)
&=&H^1_{Z\cap U}(U;\R\hom_{\A}(\A/\J,\M))\\
&=&H^1_{Z\cap U}(U;\hom_{\A}(\A/\J,\M))=0.
\eneqn
Therefore, the morphism $\phi\cl\I\to\M$ extends to a morphism $\I'\to\M$.
\qed

\proof[Proof of Lemma \ref{lem:indinj}]
Let $\{\M_i\}_{i\in I}$ be a filtrant inductive family
of injective $\A$-modules, and $\smash{\M=\indlim[i] \M_i}$.
For any coherent $\A$-module $\CF$, one has
$$\ext^k_\A(\CF,\M)\simeq\indlim[i]\ext^k_\A(\CF,\M_i),$$
and the condition (ii) in the lemma above is satisfied.
Note that any filtrant inductive limit of flabby sheaves on a Noetherian
scheme is flabby.
\qed

\begin{remark}\label{rem:injinj}
Let $\M$, $\N$ be coherent $\A$-modules,
and $Z$ a closed subset of $X$.
When $\A$ is commutative,
an injective object of
$\Mod_\qc(\A)$ is an injective object of $\Mod(\A)$, and 
$\Hom_{\DC^\rb_\qc(\A)}(\M,\R\Gamma_Z(\N)[n])$ is calculated by
$H^n(\Hom_\A(\M,\Gamma_Z(\I^\bullet)))$ for an injective
resolution $\I^\bullet$ of $\N$ in $\Mod_\qc(\A)$.
But it is not true in general.
Any injective object of $\Mod_\qc(\A)$
is a flabby sheaf (cf.\ Reamark \ref{rem:flabby}), 
and $\Gamma_Z(\I^\bullet)$ certainly calculates
$\R\Gamma_Z(\N)$. But 
$\I^\bullet$ is not necessarily a complex of
injective objects of $\Mod(\A)$, and
the functor $\Gamma_Z$
does not send injective objects of $\Mod_\qc(\A)$
to injective objects of $\Mod_\qc(\A)$.

For example, take $X=\Spec(\C[x])$, $\A=\RD_X$ 
and $\M=\RD_X/\RD_X\partial_x$,
$Z=\{0\}$.
Then one has by Hilbert's Nullstellensatz
$$\Hom_\A(\M,\Gamma_Z(\N))\simeq
\indlim[n>0]\Hom_\A(\RD_X/(\RD_X\partial_x+\RD_X x^n),\N)
=0$$
for any $\N\in \Mod_\qc(\A)$, but one has
$\Hom_{\DC^\rb_\qc(\A)}(\M,\R\Gamma_Z(\M)[2])\simeq\C$.
\end{remark}

\end{document}